\documentclass[12pt]{article}
\usepackage[sc]{mathpazo}
\usepackage[T1]{fontenc}
\usepackage[utf8]{inputenc}
\usepackage{geometry}
\usepackage{titlesec}
\titleformat{\subsection}[runin]{\bf}{\thesubsection.}{3pt}{}
\usepackage{amsmath, amsfonts, amssymb, tikz-cd, bbm}
\usepackage{graphicx}
\usepackage{caption, subcaption}
\usepackage{enumitem}

\DeclareMathOperator{\diam}{diam}

\usepackage{amsthm}
\newtheorem{theorem}{Theorem}[section]
\newtheorem{lemma}[theorem]{Lemma}
\newtheorem{corollary}[theorem]{Corollary}
\newtheorem{proposition}[theorem]{Proposition}
\theoremstyle{remark}
\newtheorem{rem}{Remark}[section]
\newtheorem{exm}{Example}[section]
\theoremstyle{definition}
\newtheorem{mydef}{Definition}[section]
\usepackage{hyperref}
\definecolor{upsPurp}{RGB}{99,0,60}
\definecolor{upsRed}{RGB}{198,11,70}
\definecolor{upsCyan}{RGB}{0,148,181}
\hypersetup{
     colorlinks=true,
     linkcolor=upsRed,
     citecolor = upsPurp,      
     urlcolor= upsCyan,
     }
\newcommand{\qeda}{\hfill \ensuremath{\Box}}

\renewcommand\thesubsection{\thesection.\Alph{subsection}}

\usepackage[
backend=biber,
style=alphabetic,
sorting=nty
]{biblatex}
\usepackage{csquotes}

\numberwithin{equation}{section}

\newcommand\blfootnote[1]{%
  \begingroup
  \renewcommand\thefootnote{}\footnote{#1}%
  \addtocounter{footnote}{-1}%
  \endgroup
}

\title{{\bf Hyperbolic embedding of infinite-dimensional convex bodies}}
\author{\sc{Yusen Long}}
\date{}

\begin{document}
\maketitle
\blfootnote{\emph{Date}: May 2023.\\ The author is supported by CSC-202108170018.}

\begin{abstract}
In this article, we use the second intrinsic volume to define a metric on the space of homothety classes of Gau{\ss}ian bounded convex bodies in a separable real Hilbert space. Using kernels of hyperbolic type, we can deduce that this space is isometrically embedded into an infinite-dimensional real hyperbolic space. Applying the Malliavin calculus, it is possible to adapt integral geometry for convex bodies in infinite dimension. Moreover, we give a new formula for computing the second intrinsic volumes of convex bodies and a characterisation of the equality case of Alexandrov-Fenchel inequality in infinite dimension, offering a description of the completion for the hyperbolic embedding of Gau{\ss}ian bounded convex bodies with dimension at least two and thus answering a question asked by Debin and Fillastre \cite{debin2018hyperbolic}.
\end{abstract}


\vspace{0.25cm}

\setcounter{tocdepth}{2}
\tableofcontents

\section{Introduction}
The idea of endowing the collection of flat figures with a hyperbolic structure dates back to W. P. Thurston's work \cite{thurston1998shapes}, where he provides a complex hyperbolic description of ``triangulations'' of a 2-sphere using flat metrics, denoted by $P$. A notable generalization of this concept is the study of subspaces of \( P \) endowed with an isometric involution, as explored in \cite{bavard1992polygones}. To begin, it is important to recognize that the homothety classes of ellipses in \(\mathbb{R}^2\) naturally correspond to the real hyperbolic plane. These subspaces are isometric to spaces of homothety classes of plane convex polygons with fixed edge directions and real hyperbolic distances. This approach has been extended to higher dimension by using mixed volumes to define hyperbolic metrics on spaces of convex polytopes in \(\mathbb{R}^n\). For \( n = 3 \), some of these spaces, which are isometric to real hyperbolic polyhedra, can be isometrically embedded into \( P \) \cite{fillastre2016remark,fillastre2017shapes}. Further, Debin and Fillastre advance this methodology by using intrinsic volumes to hyperbolise the homothety classes of Euclidean convex bodies in higher dimension \cite{debin2018hyperbolic}.

Given $2\leq d<\infty$, a convex body in $\mathbb{R}^d$ is a non-void convex compact subset of $\mathbb{R}^d$. Between any two convex bodies, it is possible to define an ``area distance'' by using the intrinsic volumes of convex bodies and mimicking the definition of the distance on the Klein model of real hyperbolic spaces. The ``area distance'' then becomes a metric on the space of homothety classes of convex bodies $K$ in $\mathbb{R}^d$ with $2\leq \dim(K)\leq d$. In \cite{debin2018hyperbolic}, Debin and Fillastre show that this metric space can be isometrically embedded into a real hyperbolic space of dimension at most $\aleph_0$. Their idea is to identify convex bodies with their support functions restricted to the unit sphere $S^{d-1}$, which reside in the Sobolev space of functions defined on $S^{d-1}$, and then compare the Sobolev subspace with the Klein model. Moreover, they show that if one normalises the convex bodies $K\subset\mathbb{R}^d$ so that $\diam(K)=1$ and positions them in the way that their Steiner points lie at $0$, then the hyperbolic embedding will be homeomorphic to the space of normalised convex bodies equipped with the Hausdorff distance.

At the end \cite[\S4]{debin2018hyperbolic}, they consider the canonical isometric embedding of $\mathbb{R}^d$ into $\mathbb{R}^{d+k}$. By identifying two convex bodies $K\subset \mathbb{R}^d$ and $K'\subset \mathbb{R}^\delta$ if $K$ only differs $K'$ from a homothety in $\mathbb{R}^{\max(d,\delta)}$, they are able to send all convex bodies of finite dimension into an infinite-dimensional real hyperbolic space. But this space is not complete. Examples are the sequence of $n$-dimensional unit balls (Example \ref{exm_balls}) and the increasing sequence of non-GB rectangles (Example \ref{exm_rect_seq}). So they ask the following question: \emph{can one give a description on the completion of the hyperbolic embedding of finite-dimensional convex bodies?}

In answering this question, the present article first gives an elementary proof for embedding homothety classes of Euclidean convex bodies into a real hyperbolic space. This proof relies on the result about kernels of (real) hyperbolic type given in \cite{monod2019self}. These kernels can be viewed as the hyperbolic analogue to kernels of positive and of conditionally negative type (compare to, for example, \cite[Appendix C]{bekka2008kazhdan}). The definition of the herein involved kernel of hyperbolic type only uses the \emph{intrinsic volumes} of convex bodies.

For finite-dimensional convex bodies, the volume of $K+rB^d\subset\mathbb{R}^d$, where $B^d$ is $d$-dimensional unit ball, is polynomial in $r>0$ (Steiner formula) and the intrinsic volumes of $K$ are defined as the normalised coefficients of this polynomial. When it comes to infinite dimension, the intrinsic volume of a convex body $K$ will be defined as the supremum of intrinsic volumes of finite-dimensional convex bodies contained in $K$. So a natural extension of the aforementioned hyperbolisation process to infinite-dimension is to consider the class of infinite-dimensional convex bodies with finite intrinsic volumes.

Let $\mathcal{H}$ be a separable Hilbert space over $\mathbb{R}$. Dudley introduces the notion of \emph{Gau{\ss}ian bounded} (abbv. \emph{GB}) subsets \cite{dudley1967sizes}. This family of subsets in $\mathcal{H}$ has been profoundly studied in the context of geometric probability and also finds its applications in ergodic theory \cite{dudley1973sample, bourgain1988almost, weber1994gb}. In \cite{chevet1976processus}, Chevet first defines the $i$-th intrinsic volume $V_i$, for $i\geq 1$, of an infinite-dimensional convex compact subset $K$ in $\mathcal{H}$ (called an \emph{infinite-dimensional convex body in $\mathcal{H}$}) by the supremum of $V_i(K')$ for all finite-dimensional convex bodies $K'\subset K$ and then shows that $K\subset \mathcal{H}$ is GB if and only if its intrinsic volumes $V_i(K)$ are finite for $i\geq 1$ (see Proposition \ref{prop_1}). So the hyperbolisation process is naturally applied to the homothety classes of GB convex bodies in $\mathcal{H}$, wherein are contained the finite-dimensional ones.

Recall that a \emph{real hyperbolic space} $\mathbb{H}_\mathbb{R}^\alpha$ is defined by the hyperboloid model constructed from a Hilbert space $\mathcal{H}'$ on $\mathbb{R}$ via a Lorentzian quadratic form (see Section \ref{sec-3.A}), where $\alpha=\dim(\mathcal{H}')$ is a cardinal and is called the \emph{dimension} of $\mathbb{H}_\mathbb{R}^\alpha$. It is a Gromov-hyperbolic space. A real hyperbolic space is of \emph{infinite dimension} if $\alpha\geq \aleph_0$, and if $\alpha=\aleph_0$, one simply writes $\mathbb{H}_\mathbb{R}^\infty$ for convenience.

\begin{theorem}\label{thm_1}
Let $\mathcal{H}$ be a separable Hilbert space over $\mathbb{R}$. Let $\mathbb{K}_2$ be the collection of homothety classes of GB convex bodies in $\mathcal{H}$ with dimension at least $2$. Then there exists an embedding $\iota:\mathbb{K}_2\hookrightarrow \mathbb{H}^\infty_\mathbb{R}$ of $\mathbb{K}_2$ into the $\aleph_0$-dimensional real hyperbolic space, thus defining a metric on $\mathbb{K}_2$, and its image $\iota(\mathbb{K}_2)$ forms a convex subset. 

Moreover, the map $\iota$ extends continuously to the homothety classes of segments in $\mathcal{H}$ and sends them to the Gromov boundary $\partial\big( \iota(\mathbb{K}_2)\big)$.
\end{theorem}

To understand GB convex bodies of infinite dimension, the techniques of infinite-dimensional analysis become indispensable.

Sudakov discovers that a GB convex body $K$ can be associated to a random variable $h_K(X)$ in $L^1(\Omega)$, where $(\Omega,\mathcal{F},\mathbb{P})$ is the probability space on which the isonormal Gau{\ss}ian process is defined, and the first intrinsic volume is the expectation of this random variable \cite{sudakov1971gaussian}. It is after decades that this random variable is recognised as the support function of the convex bodies in a separable real Hilbert space $\mathcal{H}$ and notions such as the Steiner point (or \emph{``centre''}) find their generalisations in the context of infinite-dimensional GB convex bodies \cite{vitale2001intrinsic}.

To treat these random variables, we turn to Malliavin calculus, which allows us to compute the Malliavin derivatives $Dh_K(X)$ of support function of a GB convex body, an $\mathcal{H}$-valued random variable representing an extremal point in $K$ where the isonormal Gau{\ss}ian process is maximised (see Proposition \ref{prop_der_sign}), so that $K$ can be recovered by taking the closed convex hull of the essential range of $Dh_K(X)$ (see Corollary \ref{cor_der_K}). Moreover, the Steiner point of $K$ is exactly the barycenter in $K$ with respect to the pushforward probability measure induced by $\omega\mapsto Dh_K(X)(\omega)$ (see Proposition \ref{prop_choquet}).

The support function $h_K(X)$ of any GB convex body $K\subset \mathcal{H}$ lies in the Sobolev space $\mathbb{D}^{1,2}$ and in particular, if $K$ is of finite dimension, then it recovers the Sobolev space introduced in \cite{schneider2014convex, debin2018hyperbolic}. Moreover, one can also generalise the formula from \cite[pp.298]{schneider2014convex} or \cite[Théorème 3.10]{chevet1976processus}. If we set the polarisation
$$V_2(K_1,K_2):=\frac{1}{2}\left(V_2(K_1+K_2)-V_2(K_1)-V_2(K_2)\right),$$
then we have the following result:
\begin{theorem}\label{thm_form}
Let $K,K'\subset \mathcal{H}$ be GB convex bodies and $X$ be an isonormal Gau{\ss}ian process on $\mathcal{H}$. Then
\begin{align}\label{form_mal}
    V_2(K)=\pi\mathbb{E}\left[|h_K(X)|^2-\|Dh_K(X)\|_\mathcal{H}^2\right]
\end{align}
and
\begin{align}\label{form_mal2}
    V_2(K,K')=\pi\mathbb{E}\left[h_K(X)h_{K'}(X)-(Dh_K(X),Dh_{K'}(X))_\mathcal{H}\right],
\end{align}
where $h_K,h_{K'}$ are the support functions of $K$ and $K'$.
\end{theorem}

With formulae (\ref{form_mal}) and (\ref{form_mal2}) above, we are able to characterise the equality cases of Alexandrov-Fenchel inequality for the second intrinsic volumes of GB convex bodies:

\begin{theorem}\label{thm_AF_equal}
Let $K,K'\subset \mathcal{H}$ be GB convex bodies in $\mathcal{H}$ of dimension at least $2$. Then $V_2(K,K')^2= V_2(K) V_2(K')$ if and only if $K=tK'+v$ for some $t> 0$ and $v\in\mathcal{H}$.
\end{theorem}

For polytopes, Chevet also gives general formulae to compute their intrinsic volumes of higher degree using the support function $h_P(X)$ and some other quantities \cite[Proposition 3.6']{chevet1976processus}. The technique of associating an infinite-dimensional GB convex body to its support function will allow us to work on function spaces instead of geometric objects while we try to understand GB convex bodies. But the answer to the following question remains unclear: \emph{is it possible to generalise the formulae for GB convex bodies, e.g. using the notions from Malliavin calculus?}

Recall that $L^2(\Omega)$ admits an orthogonal decomposition $\overline{\bigoplus}_{n\geq 0}\mathfrak{H}_n$ by the \emph{$n$-th Wiener chaos}. It turns out that the orthogonal projections of the support function $h_K(X)$ in $\mathfrak{H}_n$'s can completely determine the size, position and shape of $K$. The support function $h_K(X)$ of $K$ is approximated by the support functions $h_P(X)$ of polytopes contained in it. In a more general context, it is the estimation of the suprema for infinite Gau{\ss}ian processes by its finite sub-processes \cite[Chapter 13]{boucheron2013concentration}.

It is worth remarking that in \cite{debin2018hyperbolic}, Debin and Fillastre also discuss the orthogonal projection of $h_K(X)$ to $\mathfrak{H}_n$. The Malliavin calculus generalises their discussion and further gives to these projections geometric significations: the projection of $h_K(X)$ to $\mathfrak{H}_0$ is the constant function of $V_1(K)/\sqrt{2 \pi}$, the projection to $\mathfrak{H}_1$ is the Steiner point of $K$ and its projection to $\overline{\bigoplus}_{n\geq 2}\mathfrak{H}_n$ stands for its shape (see Section \ref{sec-4.B}). Nevertheless, \emph{is it possible to tell the geometric signification of the projections of $h_K(X)$ to $\mathfrak{H}_n$ for each $n\geq 2$?}

Debin and Fillastre \cite{debin2018hyperbolic} show that the homothety classes of the $n$-dimensional unit balls $[B^n]$ converge to a point $O\in\mathbb{H}^\infty_\mathbb{R}$ (see also Example \ref{exm_balls}), but they do not converge to any point in $\iota(\mathbb{K}_2)$. Using the tools from Malliavin derivative, it is possible to adapt integral geometry into infinite dimension. Particularly, we are able to answer the question asked by Debin and Fillastre: the completion for the hyperbolic embedding of finite-dimensional convex bodies is but the convex hull in $\mathbb{H}^\infty_\mathbb{R}$ of homothetic GB convex bodies and the point $O$, \emph{i.e.}
\begin{theorem}\label{thm_2}
Let $\mathcal{H}$ be a separable Hilbert space over $\mathbb{R}$ and $\iota:\mathbb{K}_2\to\mathbb{H}^\infty_\mathbb{R}$ be the embedding in Theorem \ref{thm_1}. Then $\overline{\iota(\mathbb{K}_2)}=\mathrm{co}\left(\iota(\mathbb{K}_2)\cup\{O\}\right)$ is the geodesic convex hull, or equivalently each point in $\overline{\iota(\mathbb{K}_2)}$ is uniquely associated to a function $h_K(X)+c\in\mathbb{D}^{1,2}$, where $c\geq 0$ is a constant, $K\subset\mathcal{H}$ is a GB convex body with $\mathrm{Stein}(K)=0$ and $h_K(X)$ is its support function.
\end{theorem}

\subsection*{Organisation of the article.} In Section \ref{sec-2.A}, the notions of isonormal Gau{\ss}ian process and GB convex bodies are introduced. Section \ref{sec-2.B} gives the definition of intrinsic volumes and Section \ref{sec-2.C} is dedicated to some concrete examples of GB sets. In Section \ref{sec-2.D}, we introduce the Vitale distance as well as the notion of Steiner point of infinite-dimensional GB convex bodies.

Section \ref{sec-3.A} gives an introduction to infinite-dimensional hyperbolic space. In Section \ref{sec-3.B} and Section \ref{sec-3.C}, we prove Theorem \ref{thm_1} via kernels of hyperbolic type and some other results of the hyperbolic embedding. Section \ref{sec-3.D} provides some examples of Cauchy sequences in $\iota(\mathbb{K}_2)$. 

Section \ref{sec-4.A} records some basic facts of Malliavin calculus. Section \ref{sec-4.B} treats the support function of GB convex bodies in the context of Malliavin calculus. Section \ref{sec-4.C} describes the relation between a GB convex body and its support function. The proof of Theorem \ref{thm_2} is provided in Section \ref{sec-4.D}. 

\subsection*{Acknowledgement.} The author wishes to express his sincere gratitude to Bruno Duchesne, Chenmin Sun and Richard A. Vitale for their helpful discussion and insightful comments. The author also thanks Fran\c{c}ois Fillastre for inspecting the draft of this paper.
\section{Gau{\ss}ian bounded convex bodies}
In the sequel, unless otherwise indicated, the space $\mathcal{H}$ will be referred to as a separable infinite-dimensional real Hilbert space carrying an inner product $(\cdot,\cdot)_{\mathcal{H}}$. Let $(e_i)_{i\geq 1}$ be an orthonormal basis of $\mathcal{H}$.

\subsection{GB and GC sets.}\label{sec-2.A} Recall that a \emph{centred Gau{\ss}ian process on $T$} is a collection of random variables on a probability space $(\Omega, \mathcal{F},\mathbb{P})$ indexed by $t\in T$ such that for any finite subset $\{v_1,\dots,v_n\}$ of $T$, the random vector $(X_{v_1},\dots,X_{v_n})\sim \mathcal{N}(0,\Sigma)$ is Gau{\ss}ian.

For the Hilbert space $\mathcal{H}$, a centered \textbf{isonormal Gau{\ss}ian process} is a Gaussian process $(X_v)_{v\in\mathcal{H}}$ on $\mathcal{H}$ such that $\mathbb{E}[X_v]=0$, $\mathbb{E}[X_v X_u]=(v,u)_\mathcal{H}$ for every $v,u\in\mathcal{H}$ and $X_{av+bu}=aX_v+b X_u$ for any $a,b\in\mathbb{R}$ and any $v,u\in\mathcal{H}$.

The term ``isonormal'' is due to Segal and this process is introduced in \cite{segal1954abstract} where it bears the name of ``canonical normal distribution''.

Such a process can be constructed explicitly as the following. Let $(X_i)_{i\geq 1}$ be a sequence of \emph{orthogau{\ss}ian} random variables, \emph{i.e.} independent and identically distributed random variables following $\mathcal{N}(0,1)$. Then an isonormal Gau{\ss}ian process on $\mathcal{H}$ can be defined by
$$X_v(\omega):=\sum_{i=1}^\infty (v,e_i)_\mathcal{H}X_i(\omega)$$
for all $v\in\mathcal{H}$ and all $\omega\in\Omega$. 

Conversely, for any orthonormal basis $(e_i)_{i\geq 1}$ in $\mathcal{H}$, an isonormal Gau{\ss}ian process $(X_v)_{v\in\mathcal{H}}$ on $\mathcal{H}$ must be such that $(X_{e_i})_{i\geq 1}$ are orthogau{\ss}ian. So the isonormal Gau{\ss}ian process on $\mathcal{H}$ is essentially unique.

Alternatively, it is also possible to regard the isonormal Gau{\ss}ian process as a random variable $X:\Omega\to\mathbb{R}^\mathbb{N}$ and identify $\mathcal{H}$ with $\ell^2$, so that the Gau{\ss}ian variable $X_v=(v,X)=\sum_{i=1}^\infty v_i X_i\sim\mathcal{N}(0,\|v\|^2)$.

It is a consequence of Weil's converse to Haar theorem \cite[Appendice]{andre1965integration} that in infinite-dimensional real Hilbert space $\mathcal{H}$, there is no complete analogue to Lebesgue or Haar measure. But there is still a need for measuring subsets in $\mathcal{H}$. In view of this, Dudley introduced the following class of subsets that are suitable for measuring \cite{dudley1967sizes}:
\begin{mydef}[Gau{\ss}ian bounded sets]
Let $K$ be a subset of $\mathcal{H}$ and $(X_v)_{v\in\mathcal{H}}$ be an isonormal Gau{\ss}ian process over $\mathcal{H}$. Then $K$ is \textbf{Gau{\ss}ian bounded}, or \textbf{GB} for abbreviation, if for any countable (dense) subset $C\subset K$,
$$\mathbb{P}\left(\left\{\omega\in \Omega:\exists M\in(0,\infty)\text{ such that }\sup_{v\in C}X_v(\omega)<M\right\}\right)=1.$$
\end{mydef}

\begin{rem}
For a countable set $C$, the supremum $\sup_{v\in C}X_v$ clearly defines a random variable. Since the isonormal Gau{\ss}ian process $X$ is linear, if $K\subset \mathcal{H}$ is in addition convex, then for any dense subset $C\subset K$, $\sup_{v\in C}X_v$ is actually $\sup_{v\in K}X_v$. So being GB means that the sample function $X(\cdot, \omega)$ of the isonormal Gau{\ss}ian process is uniformly bounded on $K$ for almost all $\omega\in\Omega$. We also remark that for every countable index set $C$, the random variable $\sup_{v\in C}|X_v|$ is almost surely bounded if and only if it has a finite expectation \cite{landau1970supremum}. Moreover, for a GB convex body $K$, the random variable $\sup_{v\in C}X_v(\omega)$ does not depend on the choice of countable dense subsets in $K$ \cite{vitale2001intrinsic}. Thus, we will simply write this random variable as $\sup_{v\in K}X_v$ in the sequel.
\end{rem}

Let us focus on some properties of GB sets in a separable Hilbert space for a while. For two sets $A,B$ in the vector space $\mathcal{H}$, one can define the \textbf{Minkowski sum} by
$$A+B:=\left\{x+y\in\mathcal{H}: x\in A\text{ and }y\in B\right\}.$$
Moreover, subsets of $\mathcal{H}$ are also carrying the scalar multiplication
$$tA:=\left\{tx\in\mathcal{H}:x\in A\right\}$$
for any $t$, and when $t>0$, it is called a \textbf{dilation}. As usual, a \textbf{translation} is a map $A\mapsto A+p$ for some vector $p\in\mathcal{H}$. A finite combination of dilations and translations will then be called a \textbf{homothety}. It is obvious from the definition that the class of GB sets are stable under Minkowski additions and homotheties, as it is also remarked in \cite{dudley1971seminorms}.

Moreover, if $K$ is a GB set in $\mathcal{H}$, then so will be its convex hull, which is the collection of all convex combinations of points in $K$. A subset of the GB set $K$ in $\mathcal{H}$ is also GB. These follow directly from the definition. Since being GB is a closed condition, this implies that the closure of a GB set $K$ is also GB.

Finally, let us mention the following compactness result about GB set: \textit{every GB set in $\mathcal{H}$ is totally bounded and thus is relatively compact} \cite[Proposition 3.4]{dudley1967sizes}.

By taking the closed convex hull of a GB set, we are allowed to only focus on convex, compact, GB subsets in $\mathcal{H}$, which will be called \textbf{GB convex bodies} in $\mathcal{H}$. It is clear that a GB convex body in $\mathcal{H}$ cannot have a non-empty interior, otherwise it would contain an open ball and would not be totally bounded.

Another similar notion to GB sets is the following:
\begin{mydef}[Gau{\ss}ian continuous set]
Let $K$ be a subset of $\mathcal{H}$ and $(X_v)_{v\in\mathcal{H}}$ be an isonormal Gau{\ss}ian process over $\mathcal{H}$. Then $K$ is \textbf{Gau{\ss}ian continuous}, or \textbf{GC} for abbreviation, if for almost all $\omega\in\Omega$, the sample function $X(\cdot,\omega)$ is continuous on $K$.
\end{mydef}

It is clear that every compact GC set is GB. But the non-GC compacta amongst the GB sets are quite narrow.

Let $K$ be a subset of $\mathcal{H}$ that is convex and symmetric. For each $v\in \mathcal{H}$, define $\|v\|_K:=\sup\left\{|(u,v)_\mathcal{H}|:u\in K\right\}$. For any two bounded convex subsets $K,K'$ in $\mathcal{H}$, note $K\prec K'$ if $K\subset \mathrm{span}(K')$ and $K$ is compact for $\|\cdot\|_{s(K')}$, where $s(K')$ is the symmetric closed convex hull of $K'$. Then a \emph{maximal} GB set is such that $K'$ will never be GB whenever $K\prec K'$. As a result, \emph{every GB set is either maximal or GC} \cite[Theorem 4.7]{dudley1967sizes}.

\subsection{Intrinsic volumes.}\label{sec-2.B} Let $K$ be a subset of $\mathcal{H}$. The dimension of $K$ will be defined by the dimension of the subspaces in $\mathcal{H}$ generated by $K$. If $K$ is a convex body in $\mathcal{H}$ of dimension $d<\infty$, then it can be identified with a convex body in $\mathbb{R}^d$ and its \textbf{$k$-th intrinsic volume}, denoted by $V_k(K)$, is a positive function that can be characterised by the \emph{Steiner formula}
\begin{align}\label{form.steiner}
\mathrm{vol}_d\left(K+r B^d\right)=\sum_{k=0}^d r^{d-k} \kappa_{d-k} V_k(K),
\end{align}
where $\mathrm{vol}_d$ is the Lebesgue measure in $\mathbb{R}^d$, $B^d$ is the unit ball in $\mathbb{R}^d$ and $\kappa_k$ is the Lebesgue measure of the unit $k$-ball. Although Steiner formula depends on the dimension of the ambient Euclidean space, we emphasise that the intrinsic volumes \emph{are independent of the dimension of the ambient Euclidean space}. We remark that if $K$ is $d$-dimensional, then $V_d(K)$ is its Lebesgue measure in $\mathbb{R}^d$.

The polarisation of $V_2$ by setting
$$V_2(K_1,K_2):=\frac{1}{2}\left(V_2(K_1+K_2)-V_2(K_1)-V_2(K_2)\right)$$
provides a positively bilinear form on finite-dimensional convex bodies. More precisely, $V_2(\cdot,\cdot)$ enjoys the following properties for finite-dimensional convex bodies:
\begin{enumerate}[label=(A\arabic*), topsep=0pt, itemsep=-1ex, partopsep=1ex, parsep=1ex ]
	\item $V_2(K,K)=V_2(K)$.\label{A1}
    \item $V_2(K_1,K_2)=V_2(K_2,K_1)$\label{A2}
    \item For any $ t\geq 0$, $V_2( t K_1,K_2)= t V_2(K_1,K_2)$.\label{A3}
    \item $V_2(K_1+K_2,K_3)=V_2(K_1,K_3)+V_2(K_2,K_3)$.\label{A4}
    \item $K_1\subseteq K_2$ implies $V_2(K_1,K_3)\leq V_2(K_2,K_3)$.\label{A5}
    \item $V_2(K_1,K_2)\geq 0$ and the equality holds if and only if $K_1$ or $K_2$ is a point, or both are segments in the same direction.\label{A6}
    \item (\textbf{Alexandrov-Fenchel}) for convex bodies $K_1$ and $K_2$ in $\mathbb{R}^d$ of dimension at least $2$, $V_2(K_1,K_2)^2\geq V_2(K_1)V_2(K_2)$, and equality holds if and only if $K_1= t K_2 +v$ for some $ t>0$ and $v\in\mathbb{R}^d$.\label{A7}
\end{enumerate}
These properties are classical and can be shown in geometrical means, see for example \cite[\S 5]{schneider2014convex}. In terms of \emph{mixed volume} or \emph{querma{\ss}integrals}, $V_2(K,K')=V_2(K,K',B^d,\dots, B^d)$ if $K,K'\subset \mathbb{R}^d$, where $B^d$ is the unit ball in $\mathbb{R}^d$ \cite[\S 5 \& \S 6]{schneider2014convex}.

If $K$ is a convex body in $\mathcal{H}$ with infinite dimension, then following \cite{chevet1976processus}, its $k$-th intrinsic volume will be defined by
$$V_k(K):=\sup\left\{V_k(C):C\subset K \text{ convex body with }\dim(C)<\infty\right\}.$$

With this definition, the following properties give a full characterisation of GB convex bodies in $\mathcal{H}$ \cite[Proposition 4.1]{chevet1976processus}:
\begin{proposition}\label{prop_1}
Let $K$ be a convex body in $\mathcal{H}$. Then the following assertions are equivalent:
\begin{enumerate}[label=(\roman{*}), topsep=0pt, itemsep=-1ex, partopsep=1ex, parsep=1ex ]
\item $K$ is GB.
\item For all $k\geq 1$, the intrinsic volume $V_k(K)$ is finite.
\item There exists a $k\geq 1$ such that the intrinsic volume $V_k(K)$ is finite.
\end{enumerate}
\end{proposition}

\begin{rem}\label{rem_orth_gp}
The definition of the intrinsic volumes of a finite-dimensional convex body $K$ in $\mathcal{H}$ does not rely on the choice of the ambient finite-dimensional subspace. So the intrinsic volumes are invariant under the action of the orthogonal group $O(\mathcal{H})$. The same holds for infinite-dimensional GB convex bodies.
\end{rem}

In particular, the first intrinsic volume of a GB convex body $K$ is given by 
\begin{align}\label{form_sudakov}
V_1(K)=\sqrt{2\pi}\mathbb{E}\left[\sup_{v\in K}X_v\right]<\infty,
\end{align}
see \cite{sudakov1971gaussian, badrikian2006mesures}. More generally, the following formula due to Tsirelson (called \emph{Kubota-Tsirelson formula}) provides a way to define or to compute the intrinsic volumes of a GB convex body without approximating it by finite-dimensional convex bodies \cite{tsirelson1985geometric}:
\begin{align}
V_k(K)=\frac{(2\pi)^{k/2}}{ k!\kappa_k}\mathbb{E}\left[\lambda_k\left(\left\{(X^1_v,\dots,X^k_v)\in\mathbb{R}^k:v\in K\right\}\right)\right],
\end{align}
where $X^1,\cdots,X^k$ are $k$ independent isonormal Gau{\ss}ian process on $\mathcal{H}$ and $\lambda_k$ is the Lebesgue measure on $\mathbb{R}^k$. In particular, we remark that it is direct from Kubota-Tsirelson formula that $V_1$ is additive.

By passing to the limit, properties \ref{A1}-\ref{A7} also hold for all GB convex bodies (but not the equality condition in \ref{A7}).

\subsection{Examples.}\label{sec-2.C} In this section, we record three examples of infinite-dimensional GB convex bodies given in \cite{dudley1967sizes}. Again, let $(e_i)_{i\geq 1}$ be an orthonormal basis of the separable Hilbert space $\mathcal{H}$.

Let $(b_i)_{i\geq 1}$ be a sequence of non-negative real numbers. We define the associated \emph{ellipsoid} by
$$E(\big(b_i),(e_i)\big)=\left\{x=\sum_{i\geq 1}x_i e_i\in\mathcal{H}: \sum_{b_i>0} \frac{x_i^2}{b_i^2}\leq 1\right\}.$$
We remark that such an ellipsoid is compact if and only if $b_i>0$ converges to $0$. It is shown that $E(\big(b_i),(e_i)\big)$ is GB if and only if $(b_i)_{i\geq 1}\in \ell^2$, or if $E$ is a \emph{Schmidt ellipsoid} \cite[Proposition 6.3]{dudley1967sizes}. Again, the closed unit ball in $\mathcal{H}$ is not GB.

Also, for a sequence $(\ell_i)_{i\geq 1}$ of positive real numbers, we define the \emph{rectangle} by
$$R\big((\ell_i),(e_i)\big)=\left\{x=\sum_{i\geq 1}x_i e_i\in\mathcal{H}:|x_i|\leq \ell_i/2\right\}.$$
Similarly, the rectangle $R\big((\ell_i),(e_i)\big)$ is GB if and only if $(\ell_i)_{i\geq 1}\in \ell^1$ \cite[Proposition 6.6]{dudley1967sizes}. Moreover, we can explicitly compute out the intrinsic volumes of a GB rectangle:
\begin{lemma}\label{lem.rect}
Let $R=R\big((\ell_i),(e_i)\big)$ be a GB rectangle as above. Then
$$V_k(R)=\sum_{i_1<i_2<\dots<i_k}\ell_{i_1} \ell_{i_2}\cdots\ell_{i_k}.$$
\end{lemma}
\begin{proof}
Recall the formula of intrinsic volumes for orthogonal product (see \cite[(4.4.2)]{chevet1976processus} or \cite[Proposition 4.2.3, Theorem 9.7.1]{klain1997introduction} for example): for every GB convex bodies $A,B\subset\mathcal{H}$ with $(A,B)_\mathcal{H}=0$,
\begin{align}\label{eq.2.1}
    V_k(A+B)=\sum_{i+j=k}V_i(A)V_j(B),
\end{align}
where $V_0(A)=V_0(B)=1$. The result follows from an induction on $k\geq 1$.
\end{proof}

It is worth noticing that both GB ellipsoids and GB rectangles are not maximal, \emph{i.e.} they are all GC sets.

Another important example are infinite-dimensional hyperoctahedra. Let $(a_i)_{i\geq 1}$ be a sequence of positive numbers. Define
$$\mathrm{Oc}\big((a_i),(e_i)\big)=\left\{v=\sum_{i\geq 1}a_i x_i e_i\in\mathcal{H}: \sum_{i\geq 1}|x_i|=1\right\}$$
to be the symmetric closed convex hull of $\{0\}\cup\{a_ie_i\in\mathcal{H}:i\geq 1\}$. Then $\mathrm{Oc}\big((a_i),(e_i)\big)$ is GB if and only if $a_i=O\big((\log i)^{-1/2}\big)$, and it is GC if and only if $a_i=o\big((\log i)^{-1/2}\big)$.

\subsection{Vitale distance.}\label{sec-2.D} Recall that the \textbf{support function} of a convex body $K$ is defined by $h_K(x)=\sup_{v\in K}(v,x)_\mathcal{H}$. The support function can also be formally extended to $\mathbb{R}^\mathbb{N}$ by 
\begin{align}\label{supp_fct}
h_K(\vec{x})=\sup_{v\in K}\sum_{i\geq 1}(v,e_i)_\mathcal{H} x_i
\end{align}
for every $\vec{x}=(x_i)_{i\geq 1}\in\mathbb{R}^\mathbb{N}$.

The random variable $h_K(X)$ makes sense for an isonormal Gau{\ss}ian process $X$ and coincides with $\sup_{v\in K} X_v$, so $V_1(K)=\sqrt{2\pi}\mathbb{E}[h_K(X)]$.

Recall that for a GB convex body $K$, the \textbf{Steiner point} of $K$ is defined (formally) by
\begin{align}\label{stein_pt}
\mathrm{Stein}(K):=\mathbb{E}\left[h_K(X)X\right],
\end{align}
where $X$ is the isonormal Gau{\ss}ian process on $\mathcal{H}$ (see \cite{vitale2001intrinsic}). This definition of Steiner point is understood in the sense that it is uniquely determined by 
$$\big(\mathrm{Stein}(K),v\big)_\mathcal{H}=\mathbb{E}[h_K(X)X_v]\in\mathbb{R}$$
for every vector $v\in V$ (see Section \ref{sec-4.A} and Section \ref{sec-4.B} for discussions on the convergence of this expectation). For finite-dimensional convex bodies in $\mathcal{H}$, this definition is the same as the original definition introduced in \cite{grunbaum1967convex}. Moreover, we will later see that the Steiner point $\mathrm{Stein}(K)$ is exactly the barycenter of $K$ with respect to the probability supported on the extremal points $\mathrm{Ext}(K)$ and inherited from the isonormal Gau{\ss}ian process $(X_v)_{v\in\mathcal{H}}$ (see Proposition \ref{prop_choquet}).

In the context of a Hilbert space $\mathcal{H}$, the Hausdorff distance between two convex bodies $K,K'\subset \mathcal{H}$ is given by $d_\mathrm{Haus}(K,K')=\|h_K-h_{K'}\|_{L^\infty(B_{\mathcal{H}})}$, where $B_\mathcal{H}$ is the closed unit ball in $\mathcal{H}$. But this distance function is insufficient for describing the behaviours of GB convex bodies: the Steiner point is not continuous with respect to the Hausdorff distance \cite{vitale1985steiner} and neither are intrinsic volumes (a unit ball of radius $r$ always has infinite intrinsic volumes even when $r\to 0$, while $\{0\}$ is GB with $V_k(\{0\})=0$ for all $k\geq 1$).

Let $K,K'$ be two convex bodies in $\mathcal{H}$. A GB convex body $L\subset \mathcal{H}$ is said to be \emph{equalising} $K$ and $K'$ if $K\subset K'+L$ and $K'\subset K+L$. So mimicking the definition of the Hausdorff distance, Vitale defines in \cite{vitale2001intrinsic} the \textbf{Vitale distance} for GB convex bodies by
$$d_\mathrm{Vit}(K,K'):=\inf\left\{V_1(L):L\subset\mathcal{H}\text{ is GB equalising }K\text{ and }K'\right\}.$$
The function $d_\mathrm{Vit}(\cdot,\cdot)$ yields a distance on GB convex bodies and a pseudodistance on compact convex subsets of $\mathcal{H}$. 

On one hand, for any convex body $K,K'\subset\mathcal{H}$, $d_\mathrm{Hauss}(K,K')\leq d_\mathrm{Vit}(K,K')$: the diameter of $L$ must be less than $V_1(L)$ by monotonicity of $V_1$, thus must be contained in a Hilbert ball with radius $V_1(L)$, and if in addition $L$ equalises $K$ and $K'$, then $d_\mathrm{Hauss}(K,K')\leq V_1(L)$.

On the other hand, if $K,K'\subset \mathbb{R}^d$, then the GB set $L$ that equalises $K$ and $K'$ in $\mathbb{R}^d$ is at most a Euclidean ball, so $d_\mathrm{Vit}(K,K')\leq V_1(B^d) d_\mathrm{Haus}(K,K')$, where $B^d$ is the $d$-dimensional Euclidean unit ball.

Equipped with $d_\mathrm{Vit}$, both the space of all GB convex bodies and the space of all GC convex bodies are complete, and the completion of finite-dimensional convex bodies under $d_\mathrm{Vit}$ is GC convex bodies \cite[Theorem 5]{vitale2001intrinsic}. Also, there exists a constant $C>0$ such that $\| \mathrm{Stein}(K)-\mathrm{Stein}(K')\|_\mathcal{H}\leq C d_\mathrm{Vit}(K,K')$ \cite[Theorem 8]{vitale2001intrinsic}, and this also demonstrates that the Steiner point is well-defined for all GB convex bodies.

Readers can refer to \cite{vitale2001intrinsic} and \cite{le2008bounded} for further discussion on Vitale distance and its relation with oscillation of GB convex bodies.
\section{Embedding into hyperbolic space}
In this section, we will construct a distance function for the family of GB convex bodies in $\mathcal{H}$ and isometrically embed it into the infinite-dimensional real hyperbolic space and its boundary.

\subsection{Hyperbolic space.}\label{sec-3.A}
Let $J$ be an index set. We define the associated $\ell^2$-space by
$$\ell^2(J):=\left\{(x_j)_{j\in J}\in \mathbb{R}^J: \sum_{j\in J}|x_j|^2<\infty\right\}.$$
We insist that $J$ is not necessarily countable, and when it is uncountable, the summation $\sum_{j\in J}|x_j|^2$ makes sense when every but countably many $x_j$'s are null.

Let $B:\ell^2(J)\times \ell^2(J)\to \mathbb{R}$ be a symmetric bilinear form. Let us consider the vector space $\mathcal{L}:=\mathbb{R}\oplus\ell^2(J)$ carrying the bilinear form
$$B_0\big((x_0,x),(y_0,y)\big)=x_0y_0-B(x,y).$$
The associated quadratic form $B_0(x,x)$ is called a \emph{Lorentzian} quadratic form and the space $\mathcal{L}$ is a \emph{Minkowski space}.

Let $\mathbb{P}\mathcal{L}$ be the projective space of $\mathcal{L}$, \emph{i.e.} the quotient space of $\mathcal{L}\setminus\{0\}$ under the equivalent relation $x\sim tx$ for every $t\in\mathbb{R}\setminus \{0\}$ and every $x\in\mathcal{L}\setminus\{0\}$. The \emph{hyperboloid model} for the $|J|$-dimensional real hyperbolic space is given by
$$\mathbb{H}^J_\mathbb{R}:=\left\{[x]\in \mathbb{P}\mathcal{L}:B_0(x,x)>0\right\}$$
equipped with a distance function $d_\mathbb{H}:\mathbb{H}^J_\mathbb{R}\times \mathbb{H}^J_\mathbb{R}\to \mathbb{R}_+$ given by
$$d_\mathbb{H}([x],[y])=\cosh^{-1}\left(\frac{\left|B_0(x,y)\right|}{\sqrt{\left|B_0(x,x)\right|\cdot \left|B_0(y,y)\right|}}\right)$$
for every $[x],[y]\in\mathbb{H}^J_\mathbb{R}$.

In the following, we will use the notation $\mathbb{H}^\alpha_\mathbb{R}$ to denote an $\alpha$-dimensional real hyperbolic space, \emph{i.e.} $\alpha$ is the cardinal of $J$, and if $\alpha\geq \aleph_0$, we may simply call $\mathbb{H}^\alpha_\mathbb{R}$ an \textbf{infinite-dimensional real hyperbolic space}.

Recall that a \textbf{geodesic} in a metric space is an isometric embedding of a real interval. A metric space is called \emph{geodesic} if there is a geodesic segment connecting every two points in the space, and it is called \emph{uniquely geodesic} if such geodesic is unique.

The metric topology induced by $d_\mathbb{H}$ on the hyperbolic space $\mathbb{H}^\alpha_\mathbb{R}$ coincides with the quotient topology inherited from $\mathbb{P}\mathcal{L}$ and the metric space $\mathbb{H}^\alpha_\mathbb{R}$ is thus uniquely geodesic and complete \cite[Proposition 2.2.2]{das2017geometry}.

It is classical that the curvature of \emph{real} hyperbolic spaces is $-1$. Readers can refer to the proof of \cite[\S II.10, Theorem 10.10]{bridson2013metric}. Although it is done for finite-dimensional hyperbolic space, no dimensional argument is involved in the proof and thus the result also holds for $\mathbb{H}^\alpha_\mathbb{R}$. Moreover, by the virtue of \cite[Proposition 3.3.4]{das2017geometry}, the geodesic metric space $\mathbb{H}^\alpha_\mathbb{R}$ and all $\mathrm{CAT}(-1)$ spaces are Gromov hyperbolic in the sense that any side of a geodesic triangle is included in the $\log 2$-neighbourhood of the other two sides, or ``\emph{all geodesic triangles are $\log 2$-thin}'' in metric geometry jargon.

More generally, let $(X,d)$ be a metric space. The \emph{Gromov product} of two points $x,y\in X$ based at $o\in X$ is given by
$$\langle x,y\rangle_o:=\frac{1}{2}\left(d(x,o)+d(y,o)-d(x,y)\right).$$
The space $(X,d)$ is called \textbf{Gromov hyperbolic} if there exists a $\delta\geq 0$ such that 
\begin{align}\label{ineq.delta}
    \langle x,y\rangle_o\geq \min\left\{\langle x,z\rangle_o,\langle z,y\rangle_o\right\}-\delta
\end{align}
for all $x,y,z,o\in X$. This definition is coherent to the thinness of geodesic triangles mentioned above, if the concerned metric space is geodesic \cite[\S III.H]{bridson2013metric}. If the space is $\delta$-hyperbolic and geodesic, then we can interpret the Gromov product by the distance between the based point and the geodesic connecting the two points in the following sense:
\begin{align}\label{ineq1}
    d(o,[x,y])-\delta\leq \langle x,y\rangle_o\leq d(o,[x,y]),
\end{align}
where $[x,y]$ is any geodesic between $x$ and $y$, see for example \cite[III.H.1]{bridson2013metric}.

A sequence $(x_n)_{n\in\mathbb{N}}$ in a Gromov hyperbolic metric space $(X,d)$ is called \emph{Cauchy-Gromov} if $\langle x_n,x_m\rangle_o\to \infty$ as $n,m\to \infty$. This definition does not depend on the choice of the based point $o\in X$. The \textbf{Gromov boundary} $\partial X$ of $(X,d)$ is the equivalent classes of Cauchy-Gromov sequences in it and two sequences $(x_n)_{n\in\mathbb{N}}$ and $(y_n)_{n\in\mathbb{N}}$ are equivalent if $\langle x_n,y_m\rangle_o\to \infty$ as $n,m\to \infty$. We denote $\overline{X}=X\cup \partial X$. The Gromov product can be extended to $\partial X$ by defining for all $x\in X$ and $\xi\in\partial X$
$$\langle x,\xi\rangle_o:=\sup \liminf_{n\to \infty}\langle x,x_n\rangle_o,$$
where the supremum is taken among all Cauchy-Gromov sequences in $\xi \in \partial X$, and also by setting for all $\eta,\xi\in\partial X$
$$\langle \eta,\xi\rangle_o:=\sup \liminf_{n,m\to \infty}\langle y_m,x_n\rangle_o,$$
in a similar way. In particular, $\langle \xi,\xi\rangle_o=\infty$ for all $\xi \in X\cup \partial X$. Also, for any $\xi,\eta\in \partial X$ and any two sequences $x_n\to \xi$ and $y_m\to \eta$,
\begin{align}\label{estim_prod}
    \langle \xi,\eta\rangle_o-2\delta\leq \liminf_{n,m\to \infty}\langle x_n,y_m\rangle_o\leq \langle \xi,\eta\rangle_o.
\end{align}
Readers may refer to \cite[\S 7.2, Remarque 8]{pierre2013groupes}. Thus we can generalise inequality (\ref{ineq.delta}) to the entire $\overline{X}$ with a larger $\delta'>0$. We say that $x\in X$ converges to $\xi\in \partial X$ if $\langle x,\xi\rangle_o\to\infty$.

If $\gamma:[0,\infty)\to X$ is a geodesic ray in a Gromov hyperbolic space, then $\gamma(n)$ forms a Cauchy-Gromov sequence. Also, using the $\delta$-thinness of geodesic triangles, it is not hard to see that two geodesic rays converging to a same point on the Gromov boundary will eventually fall into the $\delta$-neighbourhood of each other; conversely, two geodesic rays converging to different points on the boundary are forcibly to have an infinite Hausdorff distance (see for example \cite[\S III.H]{bridson2013metric}).

Finally, we mention that complete $\mathrm{CAT}(-1)$ spaces $X$ are \textbf{regularly geodesic}: $\overline{X}$ is uniquely geodesic and the unique geodesics $[x_n,y_n]$ converge to the unique geodesic $[x,y]$ provided $x_n\to x$ and $y_n\to y$ for any pair $(x,y)\in \overline{X}\times \overline{X}$ \cite[Proposition 4.4.4]{das2017geometry}. In particular, the space $\mathbb{H}^\alpha_\mathbb{R}$ is regularly geodesic.

\subsection{Kernel of hyperbolic type.}\label{sec-3.B} Kernels of positive and of conditionally negative type are classical tools for the study of embeddings into spherical and Euclidean spaces respectively. A similar notion is also available for embedding into hyperbolic spaces.

\begin{mydef} Given a set $X$, a \textbf{kernel of (real) hyperbolic type} on $X$ is a function $\beta:X\times X\to \mathbb{R}$ that is symmetric, non-negative, equal to $1$ on the diagonal with
$$\sum_{i,j=1}^n c_i c_j \beta(x_i,x_j)\leq\left(\sum_{k=1}^n c_k\beta(x_k,x_0)\right)^2$$
for all $n\geq 1$, any $x_0,x_1,\dots,x_n\in X$ and any $c_1,\dots,c_n\in \mathbb{R}$.
\end{mydef}

By rearranging the terms, it is not hard to see that $\beta(\cdot,\cdot)$ is a kernel of hyperbolic type if and only if for every $z\in X$, the function
$$N(x,y):=\beta(x,z)\beta(y,z)-\beta(x,y)$$
is a kernel of positive type, \emph{i.e.} $\sum_{i,j=1}^n c_i c_j N(x_i,x_j)\geq 0$ for all $n\geq 1$, any $x_0,x_1,\dots,x_n\in X$ and any $c_1,\dots,c_n\in \mathbb{R}$.

Using the \emph{GNS construction} (see \cite[\S 3.B]{monod2019self} and \cite[Appendix C]{bekka2008kazhdan}), if a set $X$ is equipped with a kernel of hyperbolic type, then $X$ can be endowed with a distance function given by
\begin{align}\label{dist_gen}
d_\mathbb{H}(x,y)=\cosh^{-1}\big(\beta(x,y)\big),
\end{align}
and can be isometrically embedded into $\mathbb{H}^\alpha_\mathbb{R}$, where $\alpha$ is a cardinal that is at most $|X|$, the cardinality of $X$ \cite[Theorem 3.4]{monod2019self}.

Let $\mathcal{K}_2$ be the family of translation classes of GB convex bodies in $\mathcal{H}$ with $\dim(K)\geq 2$, \emph{i.e.} $K$ and $K'$ are identified in $\mathcal{K}_2$ if there exists a $p\in\mathcal{H}$ such that $K=K'+p$.

Two main inconveniences of treating $(\mathcal{K}_2,V_2)$ are that \ref{A3} only holds for positive numbers and that $K+(-K)\neq \{0\}$, say $B^n+(-B^n)=2B^n$. Due to these observations, the space $\mathcal{K}_2$ looks more like the positive cone of a vector space instead of the entire space. To rule out the difficulties, let $\widetilde{\mathcal{K}_2}$ be the \emph{real vector space} spanned by elements $\widetilde{K}$ for all $K\in \mathcal{K}_2$, with the identification $ t \widetilde{K}=\widetilde{ t K}$ for every $ t\geq 0$ and $\widetilde{K_1}+\widetilde{K_2}=\widetilde{K_1+K_2}$ for every $K_1,K_2\in\mathcal{K}_2$. Under this convention, we have formally $\widetilde{K}+(-\widetilde{K})=0$. As a result, any $v\in \widetilde{\mathcal{K}_2}$ can be decomposed into $\widetilde{K_1}-\widetilde{K_2}$ for some $K_1,K_2\in\mathcal{K}_2$. Furthermore, we can also linearly extend $V_2$ to $\widetilde{\mathcal{K}_2}$ by setting $V_2(-\widetilde{K_1},\widetilde{K_2})=-V_2(\widetilde{K_1},\widetilde{K_2})$.

Choose an $M\in \mathcal{K}_2$ and define
\begin{align}\label{kappa}
\rho_M(v,w):=V_2\left(v,\widetilde{M}\right)V_2\left(w,\widetilde{M}\right)-V_2(M)V_2(v,w)
\end{align}
for every pair $v,w\in \widetilde{\mathcal{K}_2}$.
\begin{proposition}\label{prop.3.1}
For any $M\in \mathcal{K}_2$, the bilinear form $\rho_M$ given as (\ref{kappa}) defines a positive semi-definite scalar product on the vector space $\widetilde{\mathcal{K}_2}$.
\end{proposition}
\begin{proof}
It is clear from the definition that $\rho_M$ is symmetric and bilinear. So it suffices to show the positive semi-definiteness. Since $\rho_{tM}=t^2\rho_M$ for every $t>0$, we may assume that $V_2(M)=1$. Take any $ t>0$ and $K_1,K_2\in\mathcal{K}_2$, by Alexandrov-Fenchel inequality, we have
\begin{align*}
0 &\leq V_2(K_1+ t M,K_2+ t M)^2-V_2(K_1+ t M)V_2(K_2+ t M)\\
  &= t^2\big[2V_2(K_1,K_2)+V_2(K_1,M)^2+V_2(K_2,M)^2-2V_2(K_1,M)V_2(K_2,M)\\
  &\qquad -V_2(K_1)-V_2(K_2)\big]+C_1(K_1,K_2)t+C_0(K_1,K_2),
\end{align*}
where $C_1(K_1,K_2)$ and $C_0(K_1,K_2)$ are constants depending only on $K_1$ and $K_2$. Because the above polynomial in $t$ always stays positive, its leading coefficient must be non-negative as well, \emph{i.e.}
\begin{align}\label{eq1}
\begin{split}
    & 2V_2(K_1,K_2)+V_2(K_1,M)^2+V_2(K_2,M)^2\geq \\
    &\qquad 2V_2(K_1,M)V_2(K_2,M)+V_2(K_1)+V_2(K_2).
\end{split}
\end{align}
Now taking any $v\in \widetilde{\mathcal{K}_2}$, we may write $v=\widetilde{K_1}-\widetilde{K_2}$ and 
\begin{align*}
    \rho_M(v,v)&=\rho_M\left(\widetilde{K_1}-\widetilde{K_2},\widetilde{K_1}-\widetilde{K_2}\right)\\
    &=\rho_M\left(\widetilde{K_1},\widetilde{K_1}\right)^2+\rho_M\left(\widetilde{K_2},\widetilde{K_2}\right)-2\rho_M\left(\widetilde{K_1},\widetilde{K_2}\right)\\
    &=2V_2(K_1,K_2)+V_2(K_1,M)^2+V_2(K_2,M)^2\\
    &\qquad -\big[2V_2(K_1,M)V_2(K_2,M)+V_2(K_1)+V_2(K_2)\big]\\
    &\geq 0.
\end{align*}
This completes the proof.
\end{proof}

Let us consider the projective space of $\widetilde{\mathcal{K}_2}$, denoted by $\mathbb{P}\widetilde{\mathcal{K}_2}$. 

Recall that two sets $K,K'\subset\mathcal{H}$ are homothetic if they $K'$ is the image of $K$ under a homothety, which is a finite combination of translations and dilations. By starting with $\widetilde{\mathcal{K}_2}$, we are taking the quotient by translations; by taking the projective space, we rule out the dilations. So the projective space $\mathbb{P}\widetilde{\mathcal{K}_2}$ contains all homothety class of GB convex bodies of dimension at least $2$. We denote by $\mathbb{K}_2\subset \mathbb{P}\widetilde{\mathcal{K}_2}$ the space of homothety classes of GB convex bodies of dimension at least $2$ and the elements in it by $[K]$ for some $K\in\mathcal{K}_2$.

\begin{proposition}\label{prop.3.2}
There is an embedding $\iota:\mathbb{K}_2\hookrightarrow\mathbb{H}^\alpha_{\mathbb{R}}$ for some cardinal $\alpha$.
\end{proposition}
\begin{proof}
By \cite[Proposition 3.3]{monod2019self} and the discussion above, it suffices to design a kernel of hyperbolic type for $\mathbb{K}_2$. We claim that
$$\beta([K_1],[K_2]):=\frac{V_2(K_1,K_2)}{\sqrt{V_2(K_1)V_2(K_2)}}$$
is of hyperbolic type. By Proposition \ref{prop.3.1}, this is equivalent to saying that for any $M\in \mathcal{K}_2$, the kernel
$$N_M([K_1],[K_2]):=\beta([K_1],[M])\beta([K_2],[M])-\beta([K_1],[K_2])=\frac{\rho_M\left(\widetilde{K_1},\widetilde{K_2}\right)}{V_2(M)\sqrt{V_2(K_1)V_2(K_2)}}$$
is of positive type. Since $V_2(M)>0$, we have
$$\sum_{i,j=1}^n c_i c_j N_M([K_i],[K_j])=\rho_M\left(\sum_{i=1}^n \frac{c_i\widetilde{K_i}}{\sqrt{V_2(K_i)}},\sum_{i=1}^n \frac{c_i\widetilde{K_i}}{\sqrt{V_2(K_i)}}\right)\Big/ V_2(M)\geq0,$$
for any $c_1,\dots,c_n\in\mathbb{R}$ and $K_1,\dots,K_2\in \mathcal{K}_2$. This finishes the proof.
\end{proof}

\begin{rem}\label{rem_3.1}
We notice that the orthogonal group $O(\mathcal{H})$ also acts on the homothety classes of GB convex bodies in $\mathcal{H}$ and, as mentioned in Remark \ref{rem_orth_gp}, this action preserves the intrinsic volumes, so it induces an isometric action of $O(\mathcal{H})$ on $\iota(\mathbb{K}_2)\subset\mathbb{H}^\alpha_\mathbb{R}$.
\end{rem}

Moreover, as per (\ref{dist_gen}), there is an explicit formula for computing the distance of homothety classes of two GB convex bodies, namely
\begin{align}\label{dist_hyp}
d_\mathbb{H}\big(\iota([K]),\iota([K'])\big)=\cosh^{-1}\left(\frac{V_2(K,K')}{\sqrt{V_2(K)V_2(K')}}\right).
\end{align}
It is immediate that $d_\mathbb{H}$ is continuous with respect to $d_\mathrm{Vit}$ for convex bodies with $V_2(K)=1$: if $L$ is a GB set body $K$ and $K'$, then $d_\mathrm{Vit}(K,K')\to 0$ implies that we can make $V_1(L)\to 0$, which further indicates that $V_2(L)\to 0$ as per \cite[(4.4.1)]{chevet1976processus} or (\ref{ineq_diam}), hence $V_2(K,K')\to 1$ and $d_\mathbb{H}\big(\iota([K]),\iota([K'])\big)\to 0$ after property \ref{A4} and \ref{A5}.

For the embedding granted by Proposition \ref{prop.3.2}, there is a minimal cardinal (\cite[\S3]{monod2019self}). Nevertheless, to conclude the minimal dimension $\alpha$, we still need some more information on the hyperbolic geometry of convex bodies.

\subsection{Hyperbolic geometry of GB convex bodies.}\label{sec-3.C} Debin and Fillastre show that the image of $\iota$ restricted to homothety classes of finite-dimensional convex bodies is geodesic \cite{debin2018hyperbolic}. The same result also holds for infinite-dimensional GB convex bodies:
\begin{proposition}\label{prop_geod}
Let $K_0,K_1\subset \mathcal{H}$ be two GB convex bodies with $\dim(K_0),\dim(K_1)\geq 2$. Then there is a unique geodesic in $\mathbb{H}^\alpha_\mathbb{R}$ connecting $\iota([K_0])$ and $\iota([K_1])$ given by $\iota([(1-t)K_0+tK_1])$ for $t\in [0,1]$.
\end{proposition}
\begin{proof}
Without loss of generality, we may assume that $V_2(K_0)=V_2(K_1)=1$. Let $a=V_2(K_0,K_1)$. By Alexandrov-Fenchel inequality, we have $a\geq 1$. Define $K_t=(1-t)K_0+tK_1$ for $t\in [0,1]$. Then we will have
\begin{align*}
&d_\mathbb{H}\big(\iota([K_0]),\iota([K_t])\big)+d_\mathbb{H}\big(\iota([K_1]),\iota([K_t])\big)\\
=&\cosh^{-1}\left(\frac{V_2(K_t,K_0)}{\sqrt{V_2(K_t)}}\right)+\cosh^{-1}\left(\frac{V_2(K_t,K_1)}{\sqrt{V_2(K_t)}}\right)\\
=&\cosh^{-1}\left(\frac{ta+(1-t)}{\sqrt{t^2+(1-t)^2+2t(1-t)a}}\right)+\cosh^{-1}\left(\frac{(1-t)a+t}{\sqrt{t^2+(1-t)^2+2t(1-t)a}}\right)\\
=&\cosh^{-1}\left(\frac{(1-t+t^2)a+(t-t^2)a^2}{t^2+(1-t)^2+2t(1-t)a}\right)=:\cosh^{-1}\big(\phi(t,a)\big).
\end{align*}
If $a=1$, then it implies that $[K_0]=[K_1]=[K_t]$ and the existence of a geodesic is automatic. So suppose that $a>1$. Since
$$\frac{\partial}{\partial t}\phi(t,a)=\frac{(a-1)a(2t-1)}{\big(t^2+(1-t)^2+2t(1-t)a\big)^2},$$
the function has maxima $\phi(1,a)=\phi(0,a)=a$. Also, we notice that $\phi(t,a)\geq 0$. It soon follows that 
$$d_\mathbb{H}\big(\iota([K_0]),\iota([K_t])\big)+d_\mathbb{H}\big(\iota([K_1]),\iota([K_t])\big)\leq \cosh^{-1}(a)=d_\mathbb{H}\big(\iota([K_0]),\iota([K_1])\big).$$
By triangle inequality, it forces the path $\big(\iota([K_t])\big)_{t\in[0,1]}$ to be the geodesic between $\iota([K_0])$ and $\iota([K_1])$ after a suitable parametrisation.
\end{proof}

A quick computation allows us to give an isometric parametrisation for the geodesic segments in $\iota(\mathbb{K}_2)$:
\begin{corollary}\label{cor_para}
Let $K_0,K_1\subset \mathcal{H}$ be two GB convex bodies with $\dim(K_0),\dim(K_1)\geq 2$. Then
$$\iota\left(\left[\frac{K_0}{V_2(K_0)}+\frac{K_1}{V_2(K_1)}\right]\right)\in\mathbb{H}^\alpha_\mathbb{R}$$
is the midpoint on the geodesic segment between $\iota([K_0])$ and $\iota([K_1])$.
\end{corollary}
\begin{proof}
By assuming $V_2(K_0)=V_2(K_1)=1$, we will have
$$d_\mathbb{H}\big(\iota([K_0]),\iota([K_0+K_1])\big)=d_\mathbb{H}\big(\iota([K_1]),\iota([K_0+K_1]))\big)=\cosh^{-1}\left(\frac{V_2(K_0,K_1)+1}{\sqrt{V_2(K_0+K_1)}}\right),$$
which completes the proof.
\end{proof}
Let $\mathbb{K}$ be the homothety classes of GB convex bodies with \emph{non-zero dimension} in $\mathcal{H}$, \emph{i.e.} GB convex bodies that do not reduce to a singleton. Then we have the following result:
\begin{corollary}\label{cor.3.5}
The embedding $\iota$ can be extended to $\mathbb{K}\hookrightarrow\overline{\mathbb{H}^\alpha_\mathbb{R}}$ and there is a bijection between the Gromov boundary $\partial\big(\iota (\mathbb{K})\big)$ and the projective space $\mathbb{P}\mathcal{H}$.
\end{corollary}
\begin{proof}
It suffices to show that homothety classes of segments are embedded in $\partial\mathbb{H}^\alpha_\mathbb{R}$. Take a GB convex body $K\subset\mathcal{H}$ with $\dim(K)\geq 2$ and $V_2(K)=1$. Let $P$ be any segment. Then for every $t\in [0,\infty)$, we have $\dim(K+tP)\geq 2$. It is clear that the path $\big(\iota([K+tP])\big)_{t\in[0,\infty)}$ is a geodesic ray, since any finite segment of it is geodesic by Proposition \ref{prop_geod} and
$$\lim_{t\to \infty}d_\mathbb{H}(\iota([K]),\iota([K+tP]))=\lim_{t\to \infty}\cosh^{-1}\left(\frac{tV_2(K,P)+1}{\sqrt{2tV_2(K,P)+1}}\right)=\infty.$$
In particular, the sequence $\big(\iota([K+nP])\big)_{n\in\mathbb{N}}$ is a Cauchy-Gromov sequence converging to a point $\partial \mathbb{H}^\alpha_\mathbb{R}$, denoted by $\iota([P])$.

We remark that $\iota([P])$ does not depend on the choice of $K$. Indeed, for two distinct GB convex bodies $K,K'\subset\mathcal{H}$ with $\dim(K),\dim(K')\geq 2$ and $V_2(K)=V_2(K')=1$, we have that
$$d_\mathbb{H}\big(\iota([K+tP],\iota([K'+tP]\big)=\cosh^{-1}\left(\frac{tV_2(K+K',P)+1}{\sqrt{(1+2tV_2(K,P))(1+2tV_2(K',P)}}\right)$$
is bounded uniformly in $t\in[0,\infty)$, thus both $\iota([K+tP])$ and $\iota([K'+tP])$ converge to the same point on $\partial\mathbb{H}^\alpha_\mathbb{R}$.

We shall now show that $\iota$ is injective. Given two segments $P,P'\subset \mathcal{H}$ that are not in the same direction, then we have
\begin{align*}
&d_\mathbb{H}\big(\iota([K+nP],\iota([K+mP']\big)\\
=&\cosh^{-1}\left(\frac{1+nV_2(K,P)+mV_2(K,P)+nmV_2(P,P')}{\sqrt{(1+2nV_2(K,P))(1+2mV_2(K,P'))}}\right)\to \infty
\end{align*}
as $n,m\to \infty$. This means that the geodesic between $\iota([K])$ and $\iota([P])$ is not fellow travelling with the one between $\iota([K])$ and $\iota([P'])$, which forces $\iota([P])\neq \iota([P'])$.

Finally, we remark that the homothety classes of $1$-dimensional convex bodies are in bijection with $\mathbb{P}\mathcal{H}\simeq\mathrm{Gr}(1,\mathcal{H})$ by sending a segment passing through the origin $0\in\mathcal{H}$ to the subspace generated by it. Hence $\mathbb{P}\mathcal{H}$ is in bijection with $\partial\big(\iota (\mathbb{K})\big)$ via $\iota$.
\end{proof}

\begin{rem} The fact that $\mathbb{P}\mathcal{H}=\partial\iota(\mathbb{K})$ is not trivial. The Hilbert space that we use to construct $\mathbb{H}^\infty_\mathbb{R}$ is an abstract Hilbert space $\widetilde{\mathcal{K}_2}$ obtained via GNS construction and the projective space $\mathbb{P}\widetilde{\mathcal{K}_2}$, which is in bijection with $\partial \mathbb{H}^\infty_\mathbb{R}$, does not \emph{a priori} have a relation with $\mathbb{P}\mathcal{H}$. It is still unclear, for example, whether $\iota$ restricted to the segments, $\iota:\mathbb{K}\setminus\mathbb{K}_2\to \partial\mathbb{H}^\infty_\mathbb{R}$ is surjective or not. The same remark implies that the action of $O(\mathcal{H})$ by isometries on $\iota(\mathbb{K}_2)$, as mentioned in Remark \ref{rem_3.1}, is also non-trivial.
\end{rem}

Since the space $\mathbb{H}^\alpha_\mathbb{R}$ is regularly geodesic, we can further generalise the result of Proposition \ref{prop_geod} to the boundary as follows:
\begin{corollary}\label{cor.geo}
For any distinct $[K_1],[K_2]\in \mathbb{K}$, the geodesic connecting $\iota([K_1])$ to $\iota([K_2])$ is the path $\iota([tK_1+(1-t)K_2])$ for $t\in(0,1)$ under a suitable parametrisation.\qeda
\end{corollary}

In the hyperbolic structure, any GB convex bodies can be approximated by the finite-dimensional convex bodies contained in it.

\begin{proposition}\label{prop_fd_appro}
Let $K$ be a GB convex body in $\mathcal{H}$ and $(K_n)_{n\geq 1}$ be a sequence of finite-dimensional convex bodies contained in $K$ such that $V_2(K_n)\to V_2(K)$ as $n\to \infty$. Then $\iota([K_n]$ converges to $\iota([K])$.
\end{proposition}
\begin{proof}
Indeed, the Alexandrov-Fenchel's inequality and the monotonicity of $V_2$ gives the following estimation
$$1\leq \frac{V_2(K_n+K)-V_2(K_n)-V_2(K)}{2\sqrt{V_2(K_n)V_2(K)}}\leq \frac{4V_2(K)-V_2(K_n)-V_2(K)}{2\sqrt{V_2(K_n)V_2(K)}}\to 1$$
as $n\to \infty$.
\end{proof}

Let $\mathbb{K}_{2,f}$ be the homothety classes of finite-dimensional convex bodies that do not reduce to a singleton or a segment. In view of Proposition \ref{prop_fd_appro}, if $\overline{\iota(\mathbb{K}_{2,f})}$ is the completion of $\iota(\mathbb{K}_{2,f})$ in $\mathbb{H}^\alpha_\mathbb{R}$, then $\iota(\mathbb{K}_2)\subset\overline{\iota(\mathbb{K}_{2,f})}$.

\begin{corollary}\label{cor.dim}
The space $\iota(\mathbb{K}_2)$ is separable and the minimal dimension $\alpha\leq \aleph_0$.
\end{corollary}
\begin{proof}
For any $d\geq 2$, the image in $\mathbb{H}^\alpha_\mathbb{R}$ of the homothety classes of the convex bodies in $\mathbb{R}^d$ is homeomorphic to the space of all convex bodies $K$ in $\mathbb{R}^d$ with $\mathrm{Stein(K)}=0$ and $V_2(K)=1$, endowed with Hausdorff distance \cite{debin2018hyperbolic}. So it is separable. It soon follows that $\iota(\mathbb{K}_{2,f})$ is a countable union of separable spaces and is thus separable. As a subspace of $\overline{\iota(\mathbb{K}_{2,f})}$, $\iota(\mathbb{K}_2)$ is also separable.
\end{proof}

\begin{proposition}\label{prop_sep}
The minimal dimension $\alpha$ for the embedding $\mathbb{K}\to\overline{\mathbb{H}^\alpha_\mathbb{R}}$ is $\aleph_0$.
\end{proposition}
\begin{proof}
Corollary \ref{cor.dim} proves one side, so it remains to show that $\alpha\geq \aleph_0$. Suppose \emph{ab absurdo} that $\alpha<\aleph_0$. Then $\partial \mathbb{H}^\alpha_\mathbb{R}$ is homeomorphic to $S^{\alpha-1}$. Let $(e_i)_{i\geq 1}$ be an orthonormal system and let $\sigma_i=\left\{x\in\mathcal{H}:x=te_i,\ t\in[0,1]\right\}$ be the corresponding unit segments. Suppose that the $\sigma_i$ are sent to $\eta_i\in \partial\mathbb{H}^\alpha_\mathbb{R}$ via $\iota$. Passing to a subsequence, we may assume by compactness of $S^{\alpha-1}$ that $(\eta_i)_{i\geq 1}$ converges to some $\eta$ in $\partial \mathbb{H}^\alpha_\mathbb{R}$.

For each $i\geq 1$, we choose a Cauchy-Gromov sequence $\big(\iota([K_i^{(n)}])\big)_{n\geq 1}$ along the geodesic $[\eta_i,\eta_{i+1}]$ that converges to $\eta_i$. By Corollary \ref{cor.geo}, the convex bodies $K_i^{(n)}$ are rectangles of the form $t\sigma_i+(1-t)\sigma_{i+1}$. For convenience reasons, fix a rectangle
$$R=\left\{x_1e_1+x_2e_2\in\mathcal{H}:-\frac{1}{2}\leq x_1,x_2\leq \frac{1}{2}\right\}.$$
By choosing $K_i^{(i)}$ so that $\langle \eta_i,\iota([K_i^{(i)}])\rangle_{\iota([R])}\geq 2^i$, the $\delta$-hyperbolicity of $\mathbb{H}^\alpha_\mathbb{R}$ yields that
$$\langle \iota([K_i^{(i)}]),\eta\rangle_{\iota([R])}\geq \min\left(\langle \eta_i,\eta\rangle_{\iota([R])},\langle \eta_i,\iota([K_i^{(i)}])\rangle_{\iota([R])}\right)-2\delta\to \infty$$
as $i\to \infty$. Hence $\big(\iota([K_i^{(i)}])\big)_{i\geq 1}$ is a Cauchy-Gromov sequence converging to $\eta$.

Let $t_i\in (0,1)$ be such that
$$K_i^{(i)}=t_i\sigma_i+(1-t_i)\sigma_{i+1}.$$
We may also assume \emph{a posteriori} that $t_i\to 0$ as $i\to \infty$. For any $\lambda\in [0,1]$, we pose
$$K^i_\lambda=\lambda K_i^{(i)}+(1-\lambda)K_{i+1}^{(i+1)}=\lambda t_i\sigma_i+(\lambda-\lambda t_i+t_{i+1}-\lambda t_{i+1})\sigma_{i+1}+(1-\lambda)t_{i+1}\sigma_{i+2}.$$
Applying Lemma \ref{lem.rect}, we have for sufficiently large $i$
\begin{align*}
  d_\mathbb{H}\big(\iota([R]),\iota([K_\lambda^i])\big)=&\cosh^{-1}\left(\frac{\frac{1}{2}\big(V_2(R+K^i_\lambda)-V_2(R)-V_2(K^i_\lambda)\big)}{\sqrt{V_2(R)V_2(K_\lambda^i)}}\right)\\
  \sim &\cosh^{-1}\left(\frac{1}{\sqrt{\lambda(1-\lambda)}}\right)\geq \cosh^{-1}(2),
\end{align*}
where we have used the asymptotic identification $t_i,t_{i+1}\sim 0$ as $i$ becomes large enough. Hence by (\ref{ineq1}), we have
\begin{align*}
    \langle \iota([K_i^{(i)}]), \iota([K_{i+1}^{(i+1)}])\rangle_{\iota([R])} \leq \min_{\lambda\in [0,1]}d_\mathbb{H}\big(\iota([R]),\iota([K_\lambda^i])\big)\sim \cosh^{-1}(2)<\infty,
\end{align*}
which contradicts to $\big(\iota([K_i^{(i)}])\big)_{i\geq 1}$ being Cauchy-Gromov.
\end{proof}

For convenience, we will simply denote in the sequel by $\mathbb{H}^\infty_\mathbb{R}$ the $\aleph_0$-dimensional real hyperbolic space, in which $\mathbb{K}_2$ is embedded.

Moreover, it is worth remarking that the homothety classes of polytopes are dense in $\iota(\mathbb{K}_2)$.
\begin{proposition}\label{prop_poly_dense}
The image of homothety classes of polytopes with dimension at least two is dense in $\iota(\mathbb{K}_2)\subset\mathbb{H}^\infty_\mathbb{R}$.
\end{proposition}
\begin{proof}
By the remark \cite[(3.9.1)]{chevet1976processus} (this is an important fact and will be used frequently in the sequel), for any GB convex body $K\subset\mathcal{H}$ with dimension at least $2$, $V_2(K)$ is the supremum amongst $V_2(P)$ of polytopes $P$ contained in $K$. Then it is possible to choose a sequence $(P_n)_{n\geq 1}$ of polytopes contained in $K$ such that $V_2(P_n)\to V_2(K)$ as $n\to \infty$. Then
\begin{align*}
d_\mathbb{H}\big(\iota([P_n]),\iota([K])\big)&=\cosh^{-1}\left(\frac{V_2(P_n+K)-V_2(P_n)-V_2(K)}{2\sqrt{V_2(P_n)V_2(K)}}\right) \\
&\leq \cosh^{-1}\left(\frac{V_2(2K)-V_2(P_n)-V_2(K)}{2\sqrt{V_2(P_n)V_2(K)}}\right)\\
&\to \cosh^{-1}(1)=0
\end{align*}
as $n\to \infty$.
\end{proof}

\subsection{Examples and non-examples of Cauchy sequences.}\label{sec-3.D} At the end of \cite{debin2018hyperbolic}, Debin and Fillastre ask about the completion of $\iota(\mathbb{K}_{2,f})$ inside of $\mathbb{H}^\infty_\mathbb{R}$. One pathological phenomenon observed by Debin and Fillastre is that there are Cauchy sequences in $\iota(\mathbb{K}_{2,f})$ that do not converge to any GB convex body. In this section, some more examples will be presented. In fact, these examples suggest that this is the only ill-behaved case.

As we remark, $\iota(\mathbb{K}_2)\subset\overline{\iota(\mathbb{K}_{2,f})}$. But the converse is not true, \emph{i.e.} $\iota(\mathbb{K}_2)$ is not the completion of $\iota(\mathbb{K}_{2,f})$ in $\mathbb{H}^\infty_\mathbb{R}$. Here are several examples:
\begin{exm}[Unit balls] \label{exm_balls}
Let $B^n\subset \mathrm{span}(e_1,\dots,e_n)$ be the unit ball of dimension $n$ in $\mathcal{H}$. Then by the Steiner formula (\ref{form.steiner}), one can compute that $V_1(B^n)=n\kappa_n/\kappa_{n-1}$ and $V_2(B^n)=(n-1)\pi$. Using the Steiner formula (\ref{form.steiner}), for $m\geq n$
$$\mathrm{vol}_m(B^n+rB^m)=\sum_{k=0}^m r^{m-k} \kappa_{m-k}V_k(B^n+B^m)=\sum_{k=0}^m (r-1)^{m-k} \kappa_{m-k}V_k(B^n),$$
and comparing the terms for $k=0,1,2$ while $r\to \infty$, we are able to compute $V_2(B^n+B^m)$ in terms of $\kappa_k$'s, so hence $V_2(B^n,B^m)/\sqrt{V_2(B^n)V_2(B^m)}$. By Stirling's approximation, we can deduce that $\big(\iota([B^n])\big)_{n\geq 2}$ is a Cauchy sequence in $\mathbb{H}^\infty_\mathbb{R}$. This is already known in \cite[\S 4]{debin2018hyperbolic}. We remark that $V_1(B^n)/\sqrt{2V_2(B^n)}\to 1$ as $n\to \infty$.
\end{exm}

\begin{exm}[Non GB rectangles]\label{exm_rect_seq} Let $(\ell_i)_{i\geq 1}$ be a sequence of strictly positive numbers. Suppose that $(\ell_i)_{i\geq 1}$ is not in $\ell^1$ so that $R\big((\ell_i),(e_i)\big)$ is not GB. But its $n$-dimensional sections $R_n:=\prod_{i=1}^n[-\ell_i/2,\ell_i/2]$ still define a Cauchy sequence in $\mathbb{H}^\infty_\mathbb{R}$ if, and only if,
$$\lim_{n\to \infty}\frac{\sum_{i=1}^n \ell_i^2}{\Big(\sum_{i=1}^n \ell_i\Big)^2}=0,$$
or equivalently if 
$$\lim_{n\to \infty}\frac{\sum_{i=1}^n \ell_i}{\sqrt{2\sum_{1\leq i<j\leq n}\ell_i\ell_j}}=\lim_{n\to \infty}\frac{V_1(R^n)}{\sqrt{2V_2(R^n)}}=1.$$
The proof of this claim is an asymptotic analysis exercise: it suffices to use 
$$\sum_{1\leq i<j\leq n}\ell_i\ell_j\sim \frac{1}{2}\left(\sum_{i=1}^n \ell_i\right)^2$$
to prove the necessity and for the sufficiency, we can deduce a contradiction by assuming
$$\sum_{i=1}^n \ell_i^2\sim k\left(\sum_{i=1}^n\ell_i\right)^2$$
for some $k\in(0,1]$. In particular, if $\ell_i$ is of at most polynomial growth, then the corresponding sections converge; but they diverge when $\ell_i$ is exponentially increasing. Moreover, given any two sequences $(a_i)_{i\geq 1}$ and $(b_i)_{i\geq 1}$ as in the claim above, we likewise define the $n$-dimensional sections $R_n$ and $R'_n$. By Cauchy-Schwarz inequality,
$$0\leq \frac{\sum_{i=1}^n a_i b_i}{\Big(\sum_{i=1}^n a_i\Big)\Big(\sum_{i=1}^n b_i\Big)}\leq \frac{\Big(\sum_{i=1}^n a_i^2\Big)^{1/2}\Big(\sum_{i=1}^n b_i^2\Big)^{1/2}}{\Big(\sum_{i=1}^n a_i\Big)\Big(\sum_{i=1}^n b_i\Big)}\to 0$$
as $n,m\to \infty$. Hence
$$\frac{V_2(R_n,R_n')}{\sqrt{V_2(R_n)V_2(R'_n)}}\sim \frac{\Big(\sum_{i=1}^n a_i\Big)\Big(\sum_{i=1}^n b_i\Big)-\sum_{i=1}^n a_i b_i}{\Big(\sum_{i=1}^n a_i\Big)\Big(\sum_{i=1}^n b_i\Big)}\to 1
$$
as $n\to \infty$. This shows that all convergent $n$-dimensional sections that are not converging to a GB rectangle converge to the same point in $\mathbb{H}^\infty_\mathbb{R}$.
\end{exm}

\begin{exm}[Comparison between rectangles and balls]
Let $I_n=\prod_{i=1}^n[-1,1]$ and $B^n$ be the unit ball as in Example \ref{exm_balls}. Both $(I_n)_{n\geq 2}$ and $(B^n)_{n\geq 2}$ define \emph{a priori} convergent sequences in $\mathbb{H}^\infty_{\mathbb{R}}$. Since $I_n$ and $B^n$ both lie in $\mathrm{span}(e_1,\dots,e_n)\simeq \mathbb{R}^n$, the identity (see \cite[pp.298]{schneider2014convex} for example)
\begin{align}\label{form_sphm}
V_2(K,K')=\frac{n-1}{2\kappa_{n-2}}\left((h_K,h_K')_{L^2(S^{n-1})}-\frac{1}{n-1}\langle\nabla h_K,\nabla h_{K'}\rangle_{L^2(S^{n-1})}\right),
\end{align}
where $\nabla$ is the gradient on $S^{n-1}$, yields
\begin{align}\label{eq.3.6}
V_2(I_n,B^n)=\frac{n-1}{2\kappa_{n-2}}\int_{S^{n-1}}h_{I_n}(v)\,\mathrm{d}v.
\end{align}
Notice that $h_{I_n}(v)=\sum_{k=1}^n|v_k|$ for every $v\in S^{n-1}$, where $v_k=(v,e_k)_\mathcal{H}$. Applying orthogonal decomposition and Fubini's theorem to (\ref{eq.3.6}), we get
$$V_2(I_n,B^n)=\sum_{k=1}^n\frac{(n-1)\pi \kappa_{n-3}}{\kappa_{n-2}}\int_{-1}^1 |t|\left(1-t^2\right)^{\frac{n-3}{2}}\,\mathrm{d}t=\frac{2n\pi\kappa_{n-3}}{\kappa_{n-2}}.$$
Hence by Stirling's approximation
$$d_\mathbb{H}\big(\iota([I_n]),\iota([B^n])\big)=\cosh^{-1}\left(\frac{2n\pi\kappa_{n-3}}{\kappa_{n-2}}\frac{1}{\sqrt{2\pi n(n-1)^2}}\right)\to 0$$
as $n\to \infty$. As a result, those two sequences converge to the same point in $\mathbb{H}^\infty_\mathbb{R}$.
\end{exm}

Now, let $O\in\mathbb{H}^\infty_\mathbb{R}$ be the limit of homothety classes of $n$-dimensional unit balls as in Example \ref{exm_balls}.

\begin{exm}[Other non-GB limits]
It is possible to construct some more Cauchy sequences of $\iota(\mathbb{K}_2)$ in $\mathbb{H}^\infty_\mathbb{R}$ that do not converge to the image of any homothety class of a GB convex body in $\mathcal{H}$. Let $I=\{t e_1:t\in[0,1]\}$ be the unit interval. Then for every $n\geq 1$, $K_n:=I+c_nB^{n+1}$ is a GB convex body in $\mathcal{H}$ of dimension at least $2$ for some $c_n>0$. Corollary \ref{cor.3.5} indicates that $\iota([K_n])$ is on the geodesic between $\iota([B^{n+1}])$ and $\iota([I])\in\partial\mathbb{H}^\infty_\mathbb{R}$. As $\mathbb{H}^\infty_\mathbb{R}$ is regularly geodesic, by choosing a suitable parameter $c_n>0$ for every $n\geq 1$ according to Corollary \ref{cor_para}, we can make sure it converges in $\mathbb{H}^\infty_\mathbb{R}$ to a point on the geodesic between $O$ and $\iota([I])$.
\end{exm}
\section{Malliavin calculus and intrinsic volumes}
Debin and Fillastre show that the hyperbolisation process can be realised by treating the Sobolev space via spherical harmonics \cite{debin2018hyperbolic}. When it comes to infinite dimension, spherical harmonics will no longer be available since the unit sphere in the Hilbert space $\mathcal{H}$ does not admit a Haar measure (due to the converse Haar's theorem of Weil \cite{andre1965integration}). In infinite dimension, Malliavin calculus becomes indispensable.

\subsection{Wiener-It{\^o} decomposition and Malliavin derivative.}\label{sec-4.A} Let $(\Omega,\mathcal{F},\mathbb{P})$ be a probability space and $L^2(\Omega)$ be the space of real $L^2$-functions on $\Omega$ with respect to the probability measure $\mathbb{P}$. Let $\mathcal{H}$ be a separable Hilbert space on $\mathbb{R}$ and $X$ be an $\mathcal{F}$-measurable isonormal Gau{\ss}ian process on $\mathcal{H}$.

Recall that the \emph{Hermite polynomials} are polynomials $(H_n)_{n\geq 0}$ determined by the recurrence relation $H_0=1$ and $H_n'(x)=nH_{n-1}(x)$ for all $n\geq 1$ and that $\mathbb{E}[H_n(Z)]=0$ for all $n\geq 1$, where $Z\sim \mathcal{N}(0,1)$ is a normal Gau{\ss}ian random variable with variance $1$. For each $n\geq 0$, the \textbf{$n$-th Wiener chaos} $\mathfrak{H}_n$ is defined as the closure in $L^2(\Omega)$ of the linear span of the set $\{H_n(X_v):v\in \mathcal{H}\}$, where $X$ is the concerned isonormal Gau{\ss}ian process on the Hilbert space $\mathcal{H}$. In particular, the space $\mathfrak{H}_0$ consists of all constant functions and $\mathfrak{H}_1=\{X_v:v\in \mathcal{H}\}$.

The terminology follows from the original article \cite{wiener1938homogeneous} where the construction is similar to the one introduced here, yet the equivalent definition appears decades later in \cite{segal1956tensor}. The following decomposition is considered due to \cite{ito1951multiple} and reader can refer to \cite[Theorem 1.1.1]{nualart2006malliavin} or \cite[pp.64]{wiener1933fourier} for a detailed proof:
\begin{theorem}[Wiener-It{\^o} decomposition]\label{thm_wi}
Let $(\Omega,\mathcal{F},\mathbb{P})$, $L^2(\Omega)$ and $\mathfrak{H}_n$ be as above. Then one has an orthogonal decomposition $L^2(\Omega)=\overline{\bigoplus}_{n\geq 0}\mathfrak{H}_n$.
\end{theorem}
\begin{rem}
Here the symbol $\overline{\bigoplus}$ refers to the Hilbert space direct sum, which is the closure of algebraic direct sum inside of a Hilbert space.
\end{rem}

Let $d>0$ be an integer and $C^\infty_P(\mathbb{R}^d)$ be the space of smooth functions on $\mathbb{R}^d$ which, together with all their partial derivatives, have at most polynomial growth. By $\mathcal{S}$ one denotes the class of random variables $\varphi:\Omega\to \mathbb{R}$ such that there exists an $n\in \mathbb{N}$, vectors $v_1, \dots,v_n\in \mathcal{H}$ and a function $f\in C^\infty_P(\mathbb{R}^d)$ verifying
$$\varphi(\omega)=f\big(X_{v_1}(\omega),\cdots,X_{v_n}(\omega)\big)$$
for almost every $\omega\in\Omega$. The random variables $\varphi\in \mathcal{S}$ are then called \textbf{smooth random variables} and the function $f\in C^\infty_P(\mathbb{R}^d)$ appearing in the definition for $\varphi\in \mathcal{S}$ is then called a \textbf{(smooth) representation of $\varphi$}.

By using the derivatives of its smooth representation, the definition of Malliavin derivative for smooth random variables soon follows:
\begin{mydef}[Malliavin derivative]
Let $\varphi(\omega)=f\big(X_{v_1}(\omega),\cdots,X_{v_n}(\omega)\big)$ be a smooth random variable as above and $f$ be its smooth representation. Then the \textbf{Malliavin derivative} $D\varphi$ of $\varphi$ is defined by the $\mathcal{H}$-valued random variable
$$D\varphi:=\sum_{j=1}^n \partial_j f\big(X_{v_1}(\omega),\cdots,X_{v_n}(\omega)\big) v_j:\Omega\to\mathcal{H},$$
where $\partial_j f$ is the $j$-th partial derivative of $f$.
\end{mydef}

Following from Cameron-Martin theorem (see for example \cite{bogachev1998gaussian}), this definition does not depend on the smooth representations of random variables.

Moreover, Malliavin derivative enjoys the integral-by-parts property, namely
\begin{align}\label{int_parts}
\mathbb{E}\left[(D\varphi,v)_\mathcal{H}\right]=\mathbb{E}\left[\varphi X_v\right]
\end{align}
for every $v\in \mathcal{H}$ and every $\varphi\in\mathcal{S}$ \cite[Lemma 1.2.1]{nualart2006malliavin}.

For each $N\geq 0$, the space $\bigoplus_{n=0}^N \mathfrak{H}_n$ consists of smooth random variables with polynomial representations of degree at most $N$. By Stone-Weierstra{\ss} Theorem, the set $\mathcal{S}\cap \bigoplus_{n=0}^N \mathfrak{H}_n$ is dense in $\bigoplus_{n=0}^N \mathfrak{H}_n$. It soon follows from Wiener-It{\^o} decomposition that $\mathcal{S}$ is dense in $L^2(\Omega)$. So it is possible to define the Malliavin derivative of an $L^2$-random variable via approximating by smooth ones.

Let $E$ and $F$ be two Banach spaces. Let $A: \mathrm{Dom}(A) \to F$ be an (unbounded) operator from $E$ to $F$, where $\mathrm{Dom}(A)$ is a subspace of $E$ on which $A$ is defined. Such an operator is called \emph{closed}, if its graph $\Gamma(A) := \{(x, Ax) \in E\times F: x \in \mathrm{Dom}(A)\}$ is closed. An operator is called \textbf{closable}, if the closure of its graph is again the graph of an operator, called the closure of $A$.

By virtue of the following result, the Malliavin derivative can be extended to the entire $L^2(\Omega)$ \cite[Proposition 1.2.1]{nualart2006malliavin}:
\begin{proposition}\label{prop_closable}
The Malliavin derivative $D:\mathcal{S}\subset L^2(\Omega)\to L^2(\Omega; \mathcal{H})$ is closable.
\end{proposition}

Abusing notation, let $D$ be again the closure of Malliavin derivative defined on $\mathcal{S}$ and by density, $D:L^2(\Omega)\to L^2(\Omega; \mathcal{H})$ gives the Malliavin derivative for all $L^2$-random variables. By passing to limits, the formula (\ref{int_parts}) also holds for every $v\in \mathcal{H}$ and every $\varphi\in L^2(\Omega)$.

The \textbf{Sobolev space} $\mathbb{D}^{1,2}$ is the Hilbert space defined by the random variables $\varphi\in L^2(\Omega)$ such that the norm
$$\|\varphi\|_{1,2}:=\left(\mathbb{E}\left[|\varphi|^2\right]+\mathbb{E}\left[\|D\varphi\|_\mathcal{H}^2\right]\right)^{1/2}$$
is finite, where the inner product is given by $\langle \psi,\varphi\rangle_{1,2}=\mathbb{E}[\psi\varphi]+\mathbb{E}\left[(D\psi,D\varphi)_\mathcal{H}\right]$. In particular, for each $n\geq 0$, $\mathfrak{H}_n\subset \mathbb{D}^{1,2}$.

Let us introduce the adjoint of the Malliavin derivative. One defines the adjoint operator $\delta$, called the \textbf{divergence operator} or \emph{Skorohod integral}, by setting taking a random variable $\delta V$ for an $\mathcal{H}$-valued random variable $V\in L^2(\Omega;\mathcal{H})$ such that
\begin{align}\label{form_div}
\mathbb{E}[X(\delta V)]=\mathbb{E}[(V,DX)_\mathcal{H}]
\end{align}
for every $X\in L^2(\Omega)$. Moreover, Meyer's inequality implies that $\delta:L^2(\Omega;\mathcal{H})\to L^2(\Omega)$ is well-defined and is a continuous operator \cite[Proposition 1.5.4]{nualart2006malliavin}.

One also defines the \textbf{Ornstein-Uhlenbeck operator} by $\bigtriangleup = \delta D:L^2(\Omega)\to L^2(\Omega)$. The following result seems to be folklore. The proof is not difficult, but it is crucial for our purpose, so a brief proof is provided:
\begin{proposition}\label{prop_eigen}
For every $n\geq 0$, the $n$-th Wiener chaos $\mathfrak{H}_n$ is the eigenspace of Ornstein-Uhlenbeck operator $\bigtriangleup$ for the eigenvalue $n$.
\end{proposition}
\begin{proof}
Let $(e_i)_{i\geq 1}$ be an orthonormal basis of $\mathcal{H}$. Let $\vec{k}$ be a multi-index on $\mathbb{N}$ and define
\begin{align}\label{Phi_k}
\Phi_{\vec{k}}:=\prod_{i\in\mathrm{supp}(\vec{k})}H_{k_i}\left(X_{e_i}\right),
\end{align}
where the $H_j$ are the Hermite polynomials and $X$ is an isonormal Gau{\ss}ian process on $\mathcal{H}$. By definition, the collection of $\Phi_{\vec{k}}$ with $|\vec{k}|=n$ is dense in $\mathfrak{H}_n$. By the recurrence relation $DH_n(X_v)=nH_{n-1}(X_v)v$, one has
$$D\Phi_{\vec{k}}=\sum_{i}k_i\Phi_{\vec{k}-\delta_i}e_i,$$
where $\delta_i(j)=\delta_{i,j}$ is the Kronecker multi-index. Recall that by \cite[(1.46)]{nualart2006malliavin}, one has $\delta(\varphi v)+(D\varphi,v)_\mathcal{H}=\varphi X_v$ for every $\varphi\in \mathcal{S}$ and every $v\in\mathcal{H}$. Hence
\begin{align*}
\bigtriangleup \Phi_{\vec{k}}&=\sum_{i}k_i\Phi_{\vec{k}-\delta_i}X_{e_i}-\sum_{i,j}k_i\big(k_j-\delta_i(j)\big)\Phi_{\vec{k}-\delta_i-\delta_j}(e_i,e_j)_\mathcal{H}\\
&=\sum_{i} k_i\left(\Phi_{\vec{k}-\delta_i}X_{e_i}-(k_i-1)\Phi_{\vec{k}-2\delta_i}\right).
\end{align*}
Since the Hermite polynomials also enjoy the relations 
$$H_{k_i-1}(x)x-(k_i-1)H_{k_i-2}(x)=H_{k_i}(x),$$
it soon follows that $\bigtriangleup \Phi_{\vec{k}}=n\Phi_{\vec{k}}$, which extends to the entire $\mathfrak{H}_n$ by density. This proves the claimed.
\end{proof}

\subsection{Support functions as random variables.}\label{sec-4.B} In conformity with Section \ref{sec-2.D}, let $K$ be a GB convex body in $\mathcal{H}$ and let $h_K$ be its support function. Then $h_K(X)$ becomes a random variable over $\Omega$. In particular, \ref{form_sudakov} implies that $K$ is a GB convex body if and only if the random variable $h_K\in L^1(\Omega)$. 

Let $P=\overline{\mathrm{co}}(v_1,\dots,v_n)$ be the closed convex hull of $n$ points $v_1,\dots,v_n\in\mathcal{H}$, \emph{i.e.} the polytope generated these points. By \cite[Théorème 3.10]{chevet1976processus},
\begin{align}\label{form_che2}
    V_2(P)=\pi \mathbb{E}\left[h_P(X)^2-\|\sigma_P(X)\|_\mathcal{H}^2\right]
\end{align}
where $\sigma_P(X)$ is the $\mathcal{H}$-valued random variable such that $\left(\sigma_P(X),X\right)=h_P(X)$, and for any two polytopes $P,P'\subset\mathcal{H}$,
\begin{align}\label{form_che}
V_2(P,P')=\pi \mathbb{E}\left[h_P(X) h_{P'}(X)-(\sigma_P(X),\sigma_{P'}(X))_\mathcal{H}\right].
\end{align}
However, since $P,P'\subset\mathbb{R}^d=\mathrm{span}(e_1,\dots,e_d)$ for some $d>0$, after \cite[\S1.7]{schneider2014convex}, it is clear that $\sigma_P(X)=\sigma_P((X_{e_1},\dots,X_{e_d}))=\nabla h_P((X_{e_1},\dots,X_{e_d}))$. As a result, (\ref{form_che}) is exactly the same as (\ref{form_sphm}).

\begin{rem}
More precisely, in \cite[\S1.7]{schneider2014convex}, it is seen that $\sigma_P(X)$ is the \emph{Fr{\'e}chet derivative} of the support function $h_K(X)$. For more information about the connection between the Malliavin derivative and the Fr{\'e}chet derivative of functions in the Wanatabe-Sobolev space, one can refer to for example \cite{clement2006duality,kruse2014strong}. But the situation here is much simpler. 
\end{rem}

Since a polytope is always bounded, $\sigma_P(X)\in L^2(\Omega;\mathcal{H})$. It soon follows from (\ref{form_che2}) that $h_P(X)\in L^2(\Omega)$. Moreover, the random variable $h_P(X)$ has a almost-everywhere differentiable representation, namely $h_K:\mathbb{R}^d\to \mathbb{R}$, it soon yields $\sigma_P(X)=Dh_K(X)$.

Although the Malliavin derivative is closable, it is in general not continuous. To generalise the formulae (\ref{form_che2}) and (\ref{form_che}), the following lemmata will be needed:
\begin{lemma}\label{lem_wk}
Let $\varphi_n\in \mathbb{D}^{1,2}$ be a sequence of random variables such that $\varphi_n\to \varphi$ in $L^2(\Omega)$. Suppose that $\sup_{n\in\mathbb{N}}\mathbb{E}[\|D\varphi_n\|_\mathcal{H}^2]<\infty$. Then $\varphi\in\mathbb{D}^{1,2}$ and $D\varphi_n$ converges weakly in $L^2(\Omega;\mathcal{H})$ to $D\varphi$, or equivalently, $\varphi_n$ converges weakly to $\varphi$ in $\mathbb{D}^{1,2}$.
\end{lemma}
\begin{proof}
Since $\varphi_n$ converges to $\varphi$ in $L^2(\Omega)$ and $\sup_{n\in\mathbb{N}}\mathbb{E}[\|D\varphi_n\|_\mathcal{H}^2]<\infty$, $(\varphi_n)_{n\in\mathbb{N}}$ is a bounded sequence in $\mathbb{D}^{1,2}$. As $\mathbb{D}^{1,2}$ is a Hilbert space, by Banach-Alaoglu theorem, one can extract a subsequence $(\varphi_{n_k})_{k\in\mathbb{N}}$ that converges weakly to a function $\phi\in \mathbb{D}^{1,2}$. In particular, the sequence $\varphi_{n_k}$ converges to $\phi$ in $L^2(\Omega)$, which indicates that $\phi=\varphi$ and $\varphi\in\mathbb{D}^{1,2}$. These arguments can be applied to any subsequence of $\varphi_n$. It follows that $\varphi_n$ converges weakly to $\varphi$ in $\mathbb{D}^{1,2}$. This proves the claimed result.
\end{proof}

\begin{lemma}\label{lem_saks}
Let $K\subset\mathcal{H}$ be a GB convex body and let $h_K(X)\in \mathbb{D}^{1,2}$ be its support function. If there exists a sequence of polytopes $(P_m)_{m\geq 1}$ such that $h_{P_m}(X)$ converges weakly to $h_K(X)$ in $\mathbb{D}^{1,2}$, then there exists another $(\widetilde{P}_m)_{m\geq 1}$ such that $h_{\widetilde{P}_m}(X)\to h_K(X)$ in $\mathbb{D}^{1,2}$.
\end{lemma}
\begin{proof}
By Banach-Saks theorem, one can extract a subsequence $P_{n_k}$ such that the Ces{\`a}ro sums
$$\widetilde{P}_m=\frac{P_{n_1}+P_{n_2}+\cdots+P_{n_m}}{m}$$
have support functions $h_{\widetilde{P}_m}(X)$ converging to $h_K(X)$ in $\mathbb{D}^{1,2}$. 
\end{proof}

\begin{proposition}\label{prop_D12}
Let $K\subset \mathcal{H}$ be GB convex bodies and $X$ be an isonormal Gau{\ss}ian process on $\mathcal{H}$. Then there exists a sequence of polytopes $(\widetilde{P}_N)_{N\geq 1}$ such that $V_2(\widetilde{P}_N)\to V_2(K)$ and $h_{\widetilde{P}_N}(X)$ converges to $h_K(X)$ in $\mathbb{D}^{1,2}$ as $N\to \infty$.
\end{proposition}
\begin{proof}
Let $(P_n)_{n\in\mathbb{N}}$ be an increasing sequence of polytopes included in a GB convex body $K$ such that $V_i(P_n)\to V_i(K)$ as $n\to \infty$ for $i=1,2$. Since $h_{P_{n}}(X)\leq h_{P_{n+1}}(X)\leq h_K(X)$, it turns out that $h_{P_{n}}(X)$ converges to $h_K(X)$ in $L^1(\Omega)$. By passing to a subsequence, we may assume that $h_{P_{n}}(X)$ converges to $h_K(X)$ almost surely. By monotone convergence theorem, $h_{P_{n}}(X)$ also converges to $h_K(X)$ in $L^2(\Omega)$. Since $K$ is bounded, there exists an $R>0$ such that the random variable taking values amongst the extremal points of $P_n$ satisfies $\|Dh_{P_n}(X)\|_\mathcal{H}<R$ for every $n\in\mathbb{N}$. By Lemma \ref{lem_wk}, $h_K(X)\in \mathbb{D}^{1,2}$ and $h_{P_n}(X)$ converges weakly to $h_K(X)$ in $\mathbb{D}^{1,2}$. By Lemma \ref{lem_saks}, it is possible to construct polytopes $\widetilde{P}_N$ contained in $K$ with $V_1(\widetilde{P}_N)\to V_1(K)$ as $N\to \infty$ and that $h_{\widetilde{P}_N}(X)$ converges to $h_K(X)$ in $\mathbb{D}^{1,2}$. Note that Proposition \ref{prop_poly_dense} asserts that $V_2(P_i,P_j)\to V_2(K)$ as $i,j\to\infty$. Moreover, by \ref{A1} and \ref{A4}, we have
$$V_2(\widetilde{P}_N)=\frac{\sum_{i,j=1}^N V_2(P_{n_i},P_{n_j})}{N^2}\to V_2(K)$$
as $N\to\infty$.
\end{proof}
\begin{rem}
In \cite[Proposition 3.6']{chevet1976processus}, a more general form for $V_k(K)$ with $k\geq 1$ is given. But it is yet unclear whether $V_k(K)$ can also be rewritten in a similar way as (\ref{form_mal}).
\end{rem}

\begin{corollary}\label{cor_der}
Let $K$ be a GB convex body in $\mathcal{H}$ and $X$ be an isonormal Gau{\ss}ian process over $\mathcal{H}$. Then there exists $M>0$ such that the Malliavin derivative of its support function satisfies $\|Dh_{K}(X)\|^2_\mathcal{H}\leq M$ almost surely.
\end{corollary}
\begin{proof}
Since $K$ is compact, there exists $M>0$ such that $K\subset B_\mathcal{H}(0,M)$. Let $\widetilde{P}_N$ be as in the proof of Proposition \ref{prop_D12}. Then $Dh_{\widetilde{P}_N}(X)$ converges to $Dh_K(X)$ in $L^2(\Omega;\mathcal{H})$. By passing to a subsequence, the convergence is with probability $1$. But $\widetilde{P}_N\subset K\subset B_\mathcal{H}(0,M)$, so $\|Dh_{\widetilde{P}_N}(X)\|^2_\mathcal{H}\leq M$. The desired result follows from letting $N$ tend to $\infty$.
\end{proof}

\begin{proof}[Proof of Theorem \ref{thm_form}]
Formulae (\ref{form_mal}) and (\ref{form_mal2}) follows directly from (\ref{form_che}) and (\ref{form_che2}) by passing to the limits.
\end{proof}

Debin and Fillastre in \cite{debin2018hyperbolic} show that the support functions of convex bodies in $\mathbb{R}^d$ of dimension at least $2$ is embedded into a Sobolev space and if in addition that their Steiner points are positioned at $0$, then the support functions restricted to $S^{d-1}$ are $L^2$-orthogonal to the eigenspace of spherical Laplacian for the minimal positive eigenvalue. This result also holds for infinite-dimensional GB convex bodies in the context of Malliavin calculus.

\begin{lemma}\label{lem_orth}
Let $K\subset \mathcal{H}$ be a GB convex bodies such that $\mathrm{Stein}(K)=0$, $X$ be an isonormal Gau{\ss}ian process over $\mathcal{H}$ and $h_K$ be the support function of $K$. Then $h_K(X)\in\overline{\bigoplus}_{n\neq 1}\mathfrak{H}_n$.
\end{lemma}
\begin{proof}
The formula (\ref{form_mal}) shows that $h_K\in\mathbb{D}^{1,2}$. By definition (\ref{stein_pt}), $\mathbb{E}[h_K(X)X]=0$, in particular, for every $v\in \mathcal{H}$, $(\mathbb{E}[h_K(X)X],v)_\mathcal{H}=\mathbb{E}[h_K(X)X_v]=0$. But $X_v$ are exactly the elements in $\mathfrak{H}_1$. This proves the desired result.
\end{proof}

Let $J_n:L^2(\Omega)\to\mathfrak{H}_n$ be the orthogonal projection onto the $n$-th Wiener chaos. Then more generally than Lemma \ref{lem_orth}, for a GB convex body $K\subset\mathcal{H}$ and an isonormal Gau{\ss}ian process $X$ over $\mathcal{H}$, the image $J_0\big(h_K(X)\big)=\mathbb{E}[h_K(X)]=V_1(K)/\sqrt{2\pi}$ and $J_1\big(h_K(X)\big)=X_{\mathrm{Stein}(K)}$.

Moreover, by Proposition \ref{prop_eigen}, the Ornstein-Uhlenbeck operator $\bigtriangleup$ is commutative with $J_n$ by linearity of the operator and orthogonality of $\mathfrak{H}_n$. At this point, it is possible to deduce the following Rayleigh's eigenvalue theorem for random variables in $\overline{\bigoplus}_{n\neq 1}\mathfrak{H}_n$:
\begin{proposition}\label{prop_ray}
Let $\varphi\in \overline{\bigoplus}_{n\neq 1}\mathfrak{H}_n$. Then
$\mathbb{E}[\|D\varphi\|_\mathcal{H}^2]\geq 2\|\varphi\|^2_{L^2(\Omega)}-2\left(\mathbb{E}[\varphi]\right)^2\geq 0$.
\end{proposition}
\begin{proof}
By orthogonality, $\varphi=\sum_{n\neq 1}J_n\varphi$. Using (\ref{form_div}), one can compute
\begin{align*}
\mathbb{E}[\|D\varphi\|_\mathcal{H}^2]+2\left(\mathbb{E}[\varphi]\right)^2&=\sum_{n\neq 1}\mathbb{E}[J_n\varphi(\bigtriangleup J_n\varphi)]+2\left(\mathbb{E}[\varphi]\right)^2\\
&=\sum_{n\neq 1}n\mathbb{E}[|J_n\varphi|^2]+2\left(J_0\varphi\right)^2\\
&=\sum_{n\geq 1}n\mathbb{E}[|J_n\varphi|^2]+2\left(J_0\varphi\right)^2\\
&\geq 2\sum_{n\neq 1}\|J_n\varphi\|^2_{L^2(\Omega)}=2\|\varphi\|^2_{L^2(\Omega)}.
\end{align*}
The inequality $2\|\varphi\|^2_{L^2(\Omega)}-2\left(\mathbb{E}[\varphi]\right)^2\geq 0$ is a direct application of Cauchy-Schwarz inequality. This completes the proof.
\end{proof}
\begin{rem}\label{remark_eigen}
Let $K\subset\mathcal{H}$ be a GB convex body with dimension at least $2$ and positioned so that $\mathrm{Stein}(K)=0$. Let $J_n:L^2(\Omega)\to\mathfrak{H}_n$ be the orthogonal projection. Then it is clear that $h_K(X)\in \overline{\bigoplus}_{n\neq 1}\mathfrak{H}_n$ and $J_0(h_K(X))=\mathbb{E}[h_K(X)]=V_1(K)/\sqrt{2\pi}>0$. Similarly to the discussion above, if $h\in \mathbb{D}^{1,2}\cap\overline{\bigoplus}_{n\geq 2}\mathfrak{H}_n$, we shall have $\mathbb{E}[\|D h\|_{\mathcal{H}}^2]\geq 2\mathbb{E}[h^2]$, hence $-V_2(\cdot,\cdot)$ defines an inner product on the space $\mathbb{D}^{1,2}\cap\overline{\bigoplus}_{n\geq 2}\mathfrak{H}_n$. This implies that the bilinear form $V_2(\cdot,\cdot)$ defined on the space of support functions of GB convex bodies with dimension at least $2$ and $\mathrm{Stein}(K)=0$ is of Lorentzian signature, from which one can construct the infinite-dimensional hyperbolic space $\mathbb{H}^\infty_\mathbb{R}$. This recovers the discussion in \cite[Proposition 2.4]{debin2018hyperbolic}, in which $\mathbb{D}^{1,2}\cap\overline{\bigoplus}_{n\geq 2}\mathfrak{H}_n$ is identified with $H^1(\mathbb{S}^{n-1})_{01}$.
\end{rem}

\subsection{Convex bodies and support functions.}\label{sec-4.C}This section mainly concerns the following question: \emph{given a function in the Sobolev space $\mathbb{D}^{1,2}$, how can one tell whether it is the support function of a GB convex body in $\mathcal{H}$?}

In finite dimension, the support function of a convex body is a function $\mathbb{R}^d\to\mathbb{R}$ that is convex, positively homogeneous and semi-lower continuous. Moreover, each of such functions uniquely define a convex body in $\mathbb{R}^d$ by
$$K=\bigcap_{x\in\mathbb{R}^d}\left\{y\in\mathbb{R}^d:x\cdot y\leq h_K(x)\right\}.$$
It turns out that $K'\subset K\subset\mathbb{R}^d$ if and only if $h_{K'}(x)\leq h_K(x)$ for all $x\in\mathbb{R}^d$.

This property can be further generalised to infinite dimension:
\begin{proposition}\label{prop_incl}
Let $K',K\subset\mathcal{H}$ be two GB convex bodies. Then $K'\subset K$ if and only if almost surely $h_{K'}(X)\leq h_K(X)$, where $X$ is an isonormal Gau{\ss}ian process on $\mathcal{H}$.
\end{proposition}
\begin{proof}
It is clear that $K'\subset K$ implies $h_{K'}(X)\leq h_K(X)$ by definition. So it remains to show the converse. Note that by taking the contraposition, the converse is equivalent to $\mathbb{P}\{h_K(X)<0\}>0$ whenever $0\notin K$. So suppose now that $K\subset\mathcal{H}$ is a GB convex body and $0\notin K$. Then there exists a unique point $v\in K$ such that $\|v\|_\mathcal{H}=\inf\left\{\|v'\|_\mathcal{H}:v'\in K\right\}$. Indeed, $v$ is the nearest point projection to $K$ and if two distinct points $v,v'\in K$ have both the least distance to $0$ amongst points in $K$, then the midpoint $v''=v/2+v'/2\in K$ would have a strictly lesser distance to $0$. By Gram-Schmidt process, let $(e_i)_{i\geq 1}$ be an orthonormal basis of $\mathcal{H}$ such that $v=t e_1$ with $t<0$. Hence 
\begin{align}\label{eq4.5}
t=\sup_{w\in K}(w,e_1)_\mathcal{H}.
\end{align}
For simplicity, we write $w_i=(w,e_i)_\mathcal{H}$. So
\begin{align}\label{eq4.6}
h_K(X)=\sup_{w\in K}\left(w_1X_{e_1}+\sum_{i=2}^\infty w_iX_{e_i}\right)\leq \sup_{w\in K}w_1 X_{e_1}+\sup_{w'\in K}\sum_{i=2}^\infty w'_iX_{e_i}.
\end{align}
When $X_{e_1}>0$, one can use (\ref{eq4.5}) to rewrite the right-hand side of (\ref{eq4.6}) into
\begin{align}\label{eq4.7}
    tX_{e_1}+\sup_{w'\in K}\sum_{i=2}^\infty w'_iX_{e_i}.
\end{align}
Since $K$ is GB, the last term of (\ref{eq4.7}) $\sup_{w'\in K}\sum_{i=2}^\infty w'_iX_{e_i}$ is almost surely finite. By independence of $X_{e_i}$ for $i\geq 1$, the two terms in (\ref{eq4.7}) are also independent. By making $X_{e_1}\gg 0$ large enough, (\ref{eq4.7}) will be negative for a positive probability, and so will be $h_K(X)$.
\end{proof}

So it immediately follows from Proposition \ref{prop_incl} that for any index set $I$, we have
\begin{align}\label{form_union}
\sup_{i\in I} h_{K_i}(X)=h_{\widetilde{K}}(X),
\end{align}
where $\widetilde{K}=\overline{\mathrm{co}}\Big(\bigcup_{i\in I}K_i\Big)$ is the closed convex hull of the union of $K_i$'s.

Let $\phi\in\mathbb{D}^{1,2}$. Define
$$\mathcal{F}_\phi :=\left\{K\subset\mathcal{H}:K\text{ is a convex body with }h_K(X)\leq \phi\right\}.$$
Then the following result can be deduced immediately from Proposition \ref{prop_incl}:
\begin{corollary}\label{cor_conv}
Let $K\subset\mathcal{H}$ be a GB convex body. Then
\begin{align}\label{form_conv}
K=\overline{\mathrm{co}}\Big(\bigcup_{P\in \mathcal{F}_{h_K(X)}}P\Big).
\end{align}
\end{corollary}
\begin{proof}
Let $K'$ be the right-hand side of (\ref{form_conv}). It is clear that $K\subset K'$ since $K\in \mathcal{F}_{h_K(X)}$. Conversely, by (\ref{form_union}),
$$h_{K'}(X)=\sup_{P\in \mathcal{F}_{h_K(X)}} h_P(X)\leq h_K(X),$$
which implies $K'\subset K$ by Proposition \ref{prop_incl}. Hence the equality is obtained.
\end{proof}
\begin{rem}
It is somehow tautological but Corollary \ref{cor_conv} also means that $K$ is the maximal convex body so that $h_K(X)$ is bounded by $\phi=h_K(X)\in \mathbb{D}^{1,2}$.
\end{rem}

For GC convex body, the isonormal Gau{\ss}ian process has almost surely a continuous sample function, so the maxima of the sample function are well-defined and attained. More generally, for GB convex body $K\subset \mathcal{H}$, we can similarly define a point $v\in K$ to be the \textbf{maximal point of $K$} at state $\omega\in\Omega$ if for every $\varepsilon>0$,
$$\sup \left\{X_w(\omega):w\in K\cap B_\mathcal{H}(v,\varepsilon)\right\}=h_K(X)(\omega).$$
Since $\mathcal{H}$ is separable and $h_K(X)<\infty$ almost surely for any $K\subset\mathcal{H}$ GB convex body, for almost every state $\omega\in\Omega$, there exists a maximal point of $K$. Moreover, for each $v\in K$, the asymptotic error
$$\lim_{\varepsilon\to 0}\sup \left\{X_w(\omega)-X_w'(\omega):w,w'\in K\cap B_\mathcal{H}(v,\varepsilon)\right\}$$
is uniformly bounded in $L^1(\Omega)$ sense \cite{vitale2001intrinsic}.

The Malliavin derivative of the supremum of a Gau{\ss}ian process is almost surely the maximal point:
\begin{proposition}\label{prop_der_sign}
Let $K\subset \mathcal{H}$ be a GB convex body and $h_K(X)$ be its support function. Then for almost every $\omega\in \Omega$, $Dh_K(X)(\omega)$ is the unique maximal point of $K$.
\end{proposition}
\begin{proof}
Note that the proof of \cite[Lemma 2.6]{kim1990cube} is  available for any non-degenerate Gau{\ss}ian process on a compact subset of a separable Banach space, with continuous covariance kernel and having at least a maximal point for almost every $\omega\in \Omega$. It soon follows that almost surely there exists a unique maximal point for the isonormal Gau{\ss}ian process $(X_v)_{v\in K}$. The fact $Dh_K(X)$ being the unique maximal point of $K$ can be deduced via approximation arguments \cite[Lemma 3.1]{decreusefond2008hitting}.
\end{proof}

Recall that for $\mathbf{f}\in L^2(\Omega;\mathcal{H})$, one can define its \textbf{essential range} by
$$\mathrm{ess}(\mathbf{f}):=\bigcap_{\mathbf{u}=\mathbf{f}\text{ a.s.}} \overline{\mathbf{u}(\Omega)}.$$
This intersection is non-void since $\mathcal{H}$ is supposed to be separable. Note that by definition the essential range remains the same for functions that only differ on a negligible set.

\begin{corollary}\label{cor_der_K}
Let $K\subset \mathcal{H}$ be a GB convex body and $h_K(X)$ be its support function. Then the closed convex hull of the essential range of the Malliavin derivative $\overline{\mathrm{co}}\big(\mathrm{ess}(Dh_K(X))\big)=K$.
\end{corollary}
\begin{proof}
Let $K'=\overline{\mathrm{co}}\big(\mathrm{ess}(Dh_K(X))\big)$. \emph{A priori} $K'\subset K$, as $Dh_K(X)\in K$ almost surely by Proposition \ref{prop_der_sign}. Conversely, note that
$$K'=\bigcap_{n\geq 1}\left\{\left[\mathrm{co}\big(\mathrm{ess}(Dh_K(X))\big)+B_\mathcal{H}\left(0,1/n\right)\right]\cap K\right\}=:\bigcap_{n\geq 1}K_n.$$
Moreover, the support function $h_{K_n}(X)$ is a decreasing sequence convergent to $h_{K'}(X)$. But for every $n\geq 1$, by the definition of maximal points $h_{K_n}(X)\geq h_K(X)$. Letting $n\to \infty$ yields $h_{K'}(X)\geq h_K(X)$, which indicates that $K\subset K'$ after Proposition \ref{prop_incl}.
\end{proof}

A Maximal point of a GB convex body $K$, if it is well-defined, is necessarily contained in the \textbf{extremal points} $\mathrm{Ext}(K)$ of $K$. Hence we can define a measurable map $\Omega \to \mathrm{Ext}(K)$ by $\omega\mapsto Dh_K(X)(\omega)$, which will further yield a probability measure $\mu_K$ on $\mathrm{Ext}(K)$ by pushing-forward.

From a functional analysis point of view, for a GB convex body in $\mathcal{H}$, \emph{Choquet's theory} also provides a natural center called the \textbf{barycenter} when a probability measure $\mu$ is defined on $\mathrm{Ext}(K)$. The barycenter is the unique point $b\in K$ such that $\int_{\mathrm{Ext}(K)} \phi(v)\,\mathrm{d}\mu(v)=\phi(b)$ for every affine function $\phi$ defined on $\mathcal{H}$. See \cite{phelps2001lectures} for more details.

The following proposition allows us to connect the stochastic point of view to this functional analysis point of view:
\begin{proposition}\label{prop_choquet}
Let $K\subset\mathcal{H}$ be a GB convex body and let $\mu_K$ be the pushforawrd measure as above. Then $\mathrm{Stein}(K)$ is the barycenter for $\mu_K$.
\end{proposition}
\begin{proof}
By Riesz representation theorem, it suffices to show that
$$\big(\mathrm{Stein}(K),w\big)_\mathcal{H}=\int_{\mathrm{Ext}(K)} \big(v,w\big)_\mathcal{H}\,\mathrm{d}\mu_K(v)$$
for any $w\in \mathcal{H}$. Since we have $\delta(\varphi v)+(D\varphi,v)_\mathcal{H}=\varphi X_v$ for every $\varphi\in \mathbb{D}^{1,2}$ and every $v\in\mathcal{H}$ \cite[(1.46)]{nualart2006malliavin}, it follows that
\begin{align*}
    \mathbb{E}\big[\big(Dh_K(X),w\big)_\mathcal{H}\big]&= \mathbb{E}\big[h_K(X)X_w-\delta\big(h_K(X)w\big)\big]\\
    &=\big(\mathrm{Stein}(K),w\big)_\mathcal{H}-\mathbb{E}\big[\delta\big(h_K(X)w\big)\big]\\
    &=\big(\mathrm{Stein}(K),w\big)_\mathcal{H}-\mathbb{E}\big[\big(h_K(X)w,D1\big)_\mathcal{H}\big]\\
    &=\big(\mathrm{Stein}(K),w\big)_\mathcal{H},
\end{align*}
where we have applied the definition of the Steiner point (\ref{stein_pt}) and the definition of the divergence operator $\delta$. As
$$\mathbb{E}\big[\big(Dh_K(X),w\big)_\mathcal{H}\big]=\int_{\mathrm{Ext}(K)} \big(v,w\big)_\mathcal{H}\,\mathrm{d}\mu_K(v)$$
following from the definition of the pushforward measure, we thus complete the proof.
\end{proof}

Now we can show Theorem \ref{thm_AF_equal} with the descriptions above.
\begin{proof}[Proof of Theorem \ref{thm_AF_equal}]
In classical spectral theory, it is well-known that we can find an orthonormal basis \((\phi_n)_{n \geq 1}\) in \( L^2(\Omega) \) such that \(\bigtriangleup \phi_n = \lambda_n \phi_n\) for every \(n \in \mathbb{N}\). This basis can be constructed using the Gram-Schmidt process from the random variables \(\Phi_{\vec{k}}\) as described in (\ref{Phi_k}). By Proposition \ref{prop_eigen}, we may assume that $\lambda_0=0$, $\phi_0=1$ and $\lambda_n\geq 1$ for all $n>0$.

Let $K,K'\subset\mathcal{H}$ be two GB convex bodies such that $V_2(K,K')^2=V_2(K)V_2(K')$. Without loss of generality, we may assume that $\mathrm{Stein}(K)=\mathrm{Stein}(K')=0$. It suffices to show that $K=tK'$ for some $t> 0$.

By eigendecomposition, we write $h_K(X)=\sum_{n\geq 1} a_n \phi_n$ and $h_{K'}(X)=\sum_{n\geq 1} b_n \phi_n$ with $a_n,b_n\in\mathbb{R}$. By Remark \ref{remark_eigen}, we have $\lambda_n>1$ for all $n\geq 1$. Using the formulae (\ref{form_mal}) and (\ref{form_mal2}), as well as substituting $h_K(X), h_{K'}(X)$ by the eigendecompositions, we can equivalently write the assumption $V_2(K,K')^2=V_2(K)V_2(K')$ as
$$\left(a_0 b_0+\sum_{n\geq 1}(1-\lambda_n) a_n b_n\right)^2=\left(a_0^2+\sum_{n\geq 1}(1-\lambda_n) a_n^2\right)\left(b_0^2+\sum_{n\geq 1}(1-\lambda_n) b_n^2\right),$$
where we have also used the definition $\bigtriangleup=\delta D$. Rearranging this formula by grouping terms into sums of squares, we then get
\begin{align}\label{eq_SOS}
\sum_{n\geq 1}(1-\lambda_n)(a_n b_0-a_0 b_n )^2=\sum_{i<j}(1-\lambda_i)(1-\lambda_j)(a_i b_j -a_j b_i)^2.
\end{align}
Since $\lambda_n>1$ for all $n\geq 1$, the left-hand side of (\ref{eq_SOS}) is non-positive, while the right-hand side of (\ref{eq_SOS}) is non-negative. This forces each square in (\ref{eq_SOS}) to vanish, \emph{i.e.} $a_i b_j=a_j b_i$ for all $i,j\geq 0$. Consequently, we have $h_K(X)= t h_{K'}(X)$ for some $t\in\mathbb{R}$. Applying (\ref{form_mal2}) again, we then have $tV_2(K')=V_2(K,K')>0$, which shows that $t>0$. Finally, Corollary \ref{cor_der_K} allows us to deduce $K= tK'$ by taking the closed convex hull of the essential range of the Malliavin derivatives.
\end{proof}

\subsection{Completion of hyperbolic embedding.}\label{sec-4.D} Let us consider any sequence $(K_n)_{n\geq 1}$ of GB convex bodies in $\mathcal{H}$ of dimension at least $2$ such that $V_2(K_n)=1$ and $\mathrm{Stein}(K_n)=0$ for all $n\geq 1$. In this section, the criteria for $\iota([K_n])$ defines a Cauchy sequence in $\mathbb{H}^\infty_\mathbb{R}$ will be treated.

Recall the following inequality from \cite[(4.4.1)]{chevet1976processus}:
\begin{align}\label{ineq_diam}
V_2(K)\leq \frac{V_1(K)^2}{2}\leq 2\pi\mathrm{diam}(K)+V_2(K).
\end{align}
Since there is a segment of length $\mathrm{diam}(K)$ contained in $K$, the monotonicity of the intrinsic volume implies $V_1(K)\geq \mathrm{diam}(K)$. Hence we can conclude the following result:

\begin{proposition}\label{prop_wc}
Let $(K_n)_{n\geq 1}$ be a sequence of GB convex bodies of dimension at least $2$ with $V_2(K_n)=1$ and $\mathrm{Stein}(K_n)=0$ for all $n\geq 1$. If $\mathrm{diam}(K_n)\to \infty$ as $n\to \infty$, then $\iota([K_n])$ will eventually leave every bounded subset of $\mathbb{H}^\infty_\mathbb{R}$.
\end{proposition}
\begin{proof}
Let $(K_n)_{\geq 1}$ be a sequence as above. Suppose that $K$ is also a GB convex body with $\dim(K)\geq2$, $V_2(K)=1$ and $\mathrm{Stein}(K)=0$. Then one has by (\ref{ineq_diam}) and $V_1(K)\geq \mathrm{diam}(K)$ the inequality
\begin{align*}
\cosh d_\mathbb{H}\big(\iota([K]),\iota([K_n])\big) &=\frac{V_2(K_n+K)-2}{2}\\
&\geq \frac{(V_1(K)+V_1(K_n))^2-4-4\pi\mathrm{diam}(K)-4\pi\mathrm{diam}(K_n)}{4}\\
&\geq \frac{(V_1(K)+\mathrm{diam}(K_n))^2-4-4\pi\mathrm{diam}(K)-4\pi\mathrm{diam}(K_n)}{4},
\end{align*}
which diverges to $\infty$ as $n\to \infty$.
\end{proof}

Let $O\in\mathbb{H}^{\infty}_\mathbb{R}$ be the limit of the image under $\iota$ of the homothety classes of the $n$-dimensional unit balls. If $K_n$ are GB convex bodies as above with $\mathrm{diam}(K_n)\to 0$ as $n\to \infty$, then they forcibly converge to $O\in \mathbb{H}^\infty_\mathbb{R}$. This result is already suggested by Example \ref{exm_balls} and Example \ref{exm_rect_seq}.
\begin{proposition}\label{prop_O}
Let $(K_n)_{n\geq 1}$ be a sequence of GB convex bodies of dimension at least $2$ such that $V_2(K_n)=1$ and $\mathrm{Stein}(K_n)=0$ for all $n\geq 1$. Let $O\in\mathbb{H}^{\infty}_\mathbb{R}$ be the limit of the image under $\iota$ of the homothety classes of the $n$-dimensional unit balls. If $\mathrm{diam}(K_n)\to 0$ as $n\to \infty$, then $\iota([K_n])$ converges to $O$ in $d_\mathbb{H}$.
\end{proposition}
\begin{proof}
Let $(K_n)_{n\geq 1}$ be such that $\mathrm{diam}(K_n)\to \infty$ as $n\to \infty$. Taking limits on every terms in (\ref{ineq_diam}) by letting $n$ tend to $\infty$ yields $V_1(K_n)\to \sqrt{2}$. Hence
\begin{align*}
1\leq \cosh \Big(d_\mathbb{H}\big(\iota([K_n]),\iota([K_m])\big)\Big)&=\frac{V_2(K_n+K_m)-2}{2}\\
&\leq \frac{(V_1(K_n)+V_1(K_m))^2-4}{4}\to 1
\end{align*}
as $n,m\to \infty$. So $\iota([K_n])$ defines a Cauchy sequence in $\mathbb{H}^\infty_\mathbb{R}$ converging to a point. If we take the sequence $K_1, B^2/V_2(B^2),K_2,B^3/V_2(B^3),\dots$, then their images under $\iota$ will become a new Cauchy sequence, which implies that $\iota([K_n])$ converges to $O\in \mathbb{H}^\infty_\mathbb{R}$.
\end{proof}

\begin{rem}
Although this proof is geometric, the point $O$ has a very specific meaning in the Sobolev space $\mathbb{D}^{1,2}$. Under the setting above, $\iota([K_n])$ yields a Cauchy sequence in $\mathbb{H}^\infty_\mathbb{R}$ if and only if $V_2(K_n,K_m)\to 1$ as $n,m\to\infty$, which is again by (\ref{form_mal2}) equivalent to
\begin{align}\label{form_equiv}
\mathbb{E}[\|Dh_{K_n}(X)-Dh_{K_m}(X)\|_\mathcal{H}^2-|h_{K_n}(X)-h_{K_m}(X)|^2]\to 0
\end{align}
as $n,m\to \infty$. If $\mathrm{diam}(K_n)\to 0$, it turns out that $\|Dh_{K_n}(X)\|_{L^\infty(\Omega;\mathcal{H})}\to 0$. Hence in view of Proposition \ref{prop_ray} and (\ref{form_mal}), $h_{K_n}(X)$ converges in $\mathbb{D}^{1,2}$ to the constant function $1/\sqrt{\pi}$. This function corresponds to the limit $O$ of the unit balls in $\mathbb{H}^\infty_\mathbb{R}$ but it is not the support function of any GB convex body in view of Corollary \ref{cor_conv} and Corollary \ref{cor_der_K}: it only bounds the singleton $\{0\}$ but differs from the support function of $\{0\}$.
\end{rem}

Now we are able to conclude Theorem \ref{thm_2}:
\begin{proof}[Proof of Theorem \ref{thm_2}] Let $(K_n)_{n\geq 1}$ be a sequence of GB convex bodies in $\mathcal{H}$. Suppose that $\dim(K_n)\geq 2$, $\mathrm{Stein}(K_n)=0$ and $V_2(K_n)=1$, such that $\iota([K_n])$ defines a Cauchy sequence in $\mathbb{H}^\infty_\mathbb{R}$. By Proposition \ref{prop_poly_dense}, we may assume that $K_n$'s are all polytopes.

Now suppose that there exists $<A<\infty$ such that $\mathrm{diam}(K_n)\leq A$ for every $n\geq 1$. Hence $V_1(K_n)$ is bounded after (\ref{ineq_diam}). By passing to a subsequence, we may assume that $V_1(K_n)$ converges. Together with Proposition \ref{prop_ray}, one gets
\begin{align*}
&\|h_{K_n}(X)-h_{K_m}(X)\|^2_{L^2(\Omega)}-(V_1(K_n)-V_1(K_m))^2/2\pi\\
\leq&\mathbb{E}[\|Dh_{K_n}(X)-Dh_{K_m}(X)\|_\mathcal{H}^2-|h_{K_n}(X)-h_{K_m}(X)|^2].
\end{align*}
Letting $n,m\to \infty$,  (\ref{form_equiv}) forces $\|h_{K_n}(X)-h_{K_m}(X)\|^2_{L^2(\Omega)}\to0$ as $n,m\to\infty$. Hence $h_{K_n}(X)$ converges to a random variable $\phi\in L^2(\Omega)$. Moreover, by assumption and Corollary \ref{cor_der}, we have $\|Dh_{K_n}\|_\mathcal{H}\leq \mathrm{diam}(K_n)<A$ almost surely, by Lemma \ref{lem_wk}, $h_{K_n}(X)$ converges to $\phi$ weakly in $\mathbb{D}^{1,2}$.

By Lemma \ref{lem_saks}, it is possible to construct another sequence $(P_m)_{m\geq 1}$ of polytopes with $\mathrm{Stein}(P_m)=0$ and $h_{P_m}(X)\to \phi$ in $\mathbb{D}^{1,2}$ as $m\to \infty$; indeed, $P_m$ is the $m$-th Ces{\`a}ro sum of a subsequence of $(K_n)_{n\geq 1}$. Moreover, by passing to a subsequence, we may assume that $h_{P_m}(X)$ (resp. $Dh_{P_m}(X)$) converges to $\phi$ (resp. $D\phi$) almost surely.

Let $\mathcal{F}_\phi$ be the collection of GB convex bodies such that $h_K(X)\leq \phi$ almost surely. Define
$$K_\phi:=\overline{\mathrm{co}}\left(\bigcup_{P\in \mathcal{F}_\phi}P\right)$$
to be the largest GB convex body such that $h_{K_\phi}(X)\leq \phi$ and define $\psi=\phi-h_{K_\phi}(X)$.

We claim that $\overline{\mathrm{co}}\big(\mathrm{ess}(D\phi)\big)\subset K_\phi$. To this end, it suffices to show that for every $v\in \mathrm{ess}(D\phi)$, the random variable satisfies $X_v\leq \phi$ almost surely. Let $S\subset \Omega$ be a subset such that $\mathbb{P}(S)=1$ with $h_{P_m}(X)\to \phi$ and $Dh_{P_m}(X)\to D\phi$ pointwisely on $S$. Then for each $v\in \mathrm{ess}(D\phi)$, there exists $\omega'\in S$ with $Dh_{P_m}(X)(\omega')\to v$ as $m\to \infty$. Moreover,
$$\mathbb{E}\left[\left|X_v-X_{Dh_{P_m}(X)(\omega')}\right|^2\right]=\|v-Dh_{P_m}(X)(\omega')\|_\mathcal{H}^2\to 0$$
as $m\to \infty$, \emph{i.e.} the convergence is in $L^2(\Omega)$. By passing to a subsequence, there exists $S_v\subset S\subset\Omega$ with $\mathbb{P}(S_v)=1$ such that 
$X_{Dh_{P_m}(X)(\omega')}(\omega)$ converges to $X_v(\omega)$ for every $\omega\in S_v$. By taking the limit along the subsequence (and abusing the notations since the limits are the same), it soon follows that
$$X_v(\omega)=\lim_{k\to \infty}X_{Dh_{P_m}(X)(\omega')}(\omega)\leq \lim_{m\to \infty}X_{Dh_{P_m}(X)(\omega)}(\omega)=\lim_{m\to\infty}h_{P_m}(X)(\omega)=\phi(\omega)$$
for every $\omega\in S_v$. Hence $\overline{\mathrm{co}}\big(\mathrm{ess}(D\phi)\big)\subset K_\phi$.
As $h_{K_\phi}(X)+\psi=\phi$, by taking the Malliavin derivative and their closed convex hull on both sides, we have by Corollary \ref{cor_der_K} that
$$K_\phi\subset\overline{\mathrm{co}}\big(\mathrm{ess}(Dh_{K_\phi}(X))\big)\subset\overline{\mathrm{co}}\big(\mathrm{ess}(Dh_{K_\phi}(X))+\mathrm{ess}(D\psi)\big)=\overline{\mathrm{co}}\big(\mathrm{ess}(D\phi)\big)\subset K_\phi.$$
But $\overline{\mathrm{co}}\big(\mathrm{ess}(Dh_{K_\phi}(X))+\mathrm{ess}(D\psi)\big)=\overline{\mathrm{co}}\big(\mathrm{ess}(Dh_{K_\phi}(X))+\overline{\mathrm{co}}\big(\mathrm{ess}(D\psi)\big)$. As a result, $K_\phi+\overline{\mathrm{co}}\big(\mathrm{ess}(D\psi)\big)\subset K_\phi$, which forces $D\psi=0$ almost surely, \emph{i.e.} $\psi$ is almost surely a constant function. Also, because $P_m$ is the Ces\`{a}ro sum of $K_n$'s, it is clear that $V_2(P_m)$ converges to $1$ as $m\to \infty$. By convergence of $h_{P_m}(X)$ to $\phi$ in $\mathbb{D}^{1,2}$, it soon yields that
$$\pi \mathbb{E}\left[\phi^2-\|D\phi\|_\mathcal{H}^2\right]=\lim_{m\to\infty}V_2(P_m)=1.$$

Alternatively, the function $\phi$ is also the limit point of the support function of the GB convex body $K_\phi+b_n B^n$, where $b_n>0$ and $B^n$ is the $n$-dimensional unit ball, with $b_n h_{B^n}(X)\to \psi$ in $\mathbb{D}^{1,2}$. But $\iota([K_\phi+b_n B^n])$ is a point on the geodesic between $\iota([K_\phi])$ and $\iota([B^n])$ in view of Proposition \ref{prop_geod}. As $\mathbb{H}^\infty_\mathbb{R}$ is regularly geodesic, $\phi$ also represents a point on the geodesic $\iota([K_\phi])$ and $O$. In particular, one remarks that $\dim(K_\phi)< 2$ is possible. If $\dim(K_\phi)=1$, then $\iota([K_\phi])\in \partial \mathbb{H}^\infty_\mathbb{R}$; or if $\dim(K_\phi)=0$, \emph{i.e.} $K_\phi=\{0\}$, in which case $h_{K_\phi}(X)=0$ and the sequence $(K_n)_{n\geq 1}$ converges to $O\in\mathbb{H}^\infty_\mathbb{R}$.
\end{proof}

\begin{rem}
As mentioned in Remark \ref{rem_3.1}, the orthogonal group $O(\mathcal{H})$ acts on $\iota(\mathbb{K}_2)$ by isometries. This action soon extends to the completion $\overline{\iota(\mathbb{K}_2)}$ and the point $O$ is the unique $O(\mathcal{H})$-invariant point.
\end{rem}

\printbibliography

@article{segal1954abstract,
 author = {Segal, Irving E},
 journal = {American Journal of Mathematics},
 number = {3},
 pages = {721--732},
 publisher = {Johns Hopkins University Press},
 title = {Abstract Probability Spaces and a Theorem of Kolmogoroff},
 volume = {76},
 year = {1954}
}

@article{dudley1967sizes,
  title={The sizes of compact subsets of Hilbert space and continuity of Gaussian processes},
  author={Dudley, Richard M},
  journal={Journal of Functional Analysis},
  volume={1},
  number={3},
  pages={290--330},
  year={1967},
  publisher={Elsevier}
}

@article{dudley1971seminorms,
  title={On seminorms and probabilities, and abstract Wiener spaces},
  author={Dudley, Richard M and Feldman, Jacob and Le Cam, Lucien},
  journal={Annals of Mathematics},
  volume={93},
  number={2},
  pages={390--408},
  year={1971},
  publisher={JSTOR}
}

@article{chevet1976processus,
  title={Processus Gaussiens et volumes mixtes},
  author={Chevet, Simone},
  journal={Zeitschrift f{\"u}r Wahrscheinlichkeitstheorie und verwandte Gebiete},
  volume={36},
  number={1},
  pages={47--65},
  year={1976},
  publisher={Springer}
}

@book{schneider2014convex,
  title={Convex Bodies: the Brunn--Minkowski Theory},
  author={Schneider, Rolf},
  edition={2},
  year={2014},
  publisher={Cambridge University Press}
}

@book{bridson2013metric,
  title={Metric Spaces of Non-positive Curvature},
  author={Bridson, Martin R and Haefliger, Andr{\'e}},
  edition={1},
  year={2013},
  publisher={Springer-Verlag}
}

@book{das2017geometry,
  title={Geometry and Dynamics in Gromov Hyperbolic Metric Spaces: with an Emphasis on Non-proper settings},
  author={Das, Tushar and Simmons, David and Urba{\'n}ski, Mariusz},
  volume={218},
  edition={1},
  year={2017},
  publisher={American Mathematical Society}
}

@book{bekka2008kazhdan,
  title={Kazhdan's Property (T)},
  author={Bekka, Bachir and de La Harpe, Pierre and Valette, Alain},
  edition={1},
  year={2008},
  publisher={Cambridge University Press}
}

@article{monod2019self,
  title={Self-representations of the M{\"o}bius group},
  author={Monod, Nicolas and Py, Pierre},
  journal={Annales Henri Lebesgue},
  volume={2},
  pages={259--280},
  year={2019}
}

@article{debin2018hyperbolic,
  author = {Debin, Clément and Fillastre, François},
  title = {Hyperbolic geometry of shapes of convex bodies},
  journal={Groups, Geometry, and Dynamics},
  volume={16},
  number={1},
  pages={115–-140},
  year={2022}
}

@article{tsirelson1985geometric,
  title={A geometric approach to maximum likelihood estimation for an infinite-dimensional Gaussian location. II},
  author={Tsirelson, Boris Semyonovich},
  journal={Theory of Probability and its Application},
  volume={30},
  number={4},
  pages={820--828},
  year={1986}
}

@book{klain1997introduction,
  title={Introduction to Geometric Probability},
  author={Klain, Daniel A and Rota, Gian-Carlo},
  year={1997},
  edition={1},
  publisher={Cambridge University Press}
}

@book{pierre2013groupes,
  title={Sur les Groupes Hyperboliques d’apr{\`e}s Mikhael Gromov},
  author={da la Harpe, Pierre and Ghys, {\'E}tienne},
  series={Progress in Mathematics},
  year={1990},
  edition={1},
  publisher={Birkh{\"a}user}
}

@article{vitale1985steiner,
  title={The Steiner point in infinite dimensions},
  author={Vitale, Richard A},
  journal={Israel Journal of Mathematics},
  volume={52},
  number={3},
  pages={245--250},
  year={1985},
  publisher={Springer}
}

@article{thurston1998shapes,
  title={Shapes of polyhedra and triangulations of the sphere},
  author={Thurston, William P},
  journal={Geometry and Topology Monographs},
  volume={1},
  year={1998}
}

@inproceedings{weber1994gb,
  title={GB and GC sets in ergodic theory},
  author={Weber, Michel},
  booktitle={Probability in Banach Spaces, 9},
  pages={129--151},
  volume={35},
  year={1994},
  publisher={Springer-Verlag}
}

@article{bourgain1988almost,
  title={Almost sure convergence and bounded entropy},
  author={Bourgain, Jean},
  journal={Israel Journal of Mathematics},
  volume={63},
  pages={79--97},
  year={1988},
  publisher={Springer}
}

@article{dudley1973sample,
  title={Sample functions of the Gaussian process},
  author={Dudley, Richard M},
  journal={The Annals of Probability},
  volume={1},
  number={1},
  pages={66--103},
  year={1973},
  publisher={Institute of Mathematical Statistics}
}

@book{andre1965integration,
  title={L'int{\'e}gration dans les groupes topologiques et ses applications},
  author={Weil, Andr{\'e}},
  year = {1965},
  edition ={2},
  publisher = {Hermann}
}

@article{sudakov1971gaussian,
  title={Gaussian random processes and measures of solid angles in Hilbert space},
  author={Sudakov, Vladimir Nikolaevich},
  journal={Doklady Akademii Nauk},
  volume={197},
  number={1},
  pages={43--45},
  year={1971}
}

@book{badrikian2006mesures,
  title={Mesures cylindriques, espaces de Wiener et fonctions al{\'e}atoires Gaussiennes},
  author={Badrikian, Albert and Chevet, Simone},
  year={1974},
  edition={1},
  publisher={Springer-Verlag}
}

@article{vitale2001intrinsic,
  title={Intrinsic volumes and Gaussian processes},
  author={Vitale, Richard A},
  journal={Advances in Applied Probability},
  volume={33},
  number={2},
  pages={354--364},
  year={2001},
  publisher={Cambridge University Press}
}

@book{grunbaum1967convex,
  title={Convex Polytopes},
  author={Gr{\"u}nbaum, Branko},
  edition={2},
  year={2003},
  publisher={Springer-Verlag}
}

@book{bogachev1998gaussian,
  title={Gaussian Measures},
  author={Bogachev, Vladimir Igorevich},
  year={1998},
  edition={1},
  publisher={American Mathematical Society}
}

@article{wiener1938homogeneous,
  title={The homogeneous chaos},
  author={Wiener, Norbert},
  journal={American Journal of Mathematics},
  volume={60},
  number={4},
  pages={897--936},
  year={1938},
  publisher={JSTOR}
}

@article{segal1956tensor,
  title={Tensor algebras over Hilbert spaces. I},
  author={Segal, Irving E},
  journal={Transactions of the American Mathematical Society},
  volume={81},
  number={1},
  pages={106--134},
  year={1956}
}

@article{ito1951multiple,
  title={Multiple wiener integral},
  author={It{\^o}, Kiyosi},
  journal={Journal of the Mathematical Society of Japan},
  volume={3},
  number={1},
  pages={157--169},
  year={1951}
}

@book{nualart2006malliavin,
  title={The Malliavin Calculus and Related Topics},
  author={Nualart, David},
  edition={2},
  year={2006},
  publisher={Springer}
}

@book{wiener1933fourier,
  title={The Fourier Integral and Certain of its Applications},
  author={Wiener, Norbert},
  edition={1},
  year={1933},
  publisher={Cambridge University Press}
}

@article{clement2006duality,
  title={A duality approach for the weak approximation of stochastic differential equations},
  author={Clement, Emmanuelle and Kohatsu-Higa, Arturo and Lamberton, Damien},
  journal={Annals of Applied Probability},
  volume={16},
  number={3},
  pages={1124--1154},
  year={2006}
}

@book{kruse2014strong,
  title={Strong and Weak Approximation of Semilinear Stochastic Evolution Equations},
  author={Kruse, Raphael},
  year={2014},
  edition={1},
  publisher={Springer}
}

@book{boucheron2013concentration,
  title={Concentration Inequalities: A Nonasymptotic Theory of Independence},
  author={Boucheron, St{\'e}phane and Lugosi, G{\'a}bor and Massart, Pascal},
  year={2013},
  publisher={Oxford University Press}
}

@article{kim1990cube,
  title={Cube root asymptotics},
  author={Kim, Jeankyung and Pollard, David},
  journal={The Annals of Statistics},
  pages={191--219},
  year={1990},
  volume={18},
  number={1}
}

@article{decreusefond2008hitting,
  title={Hitting Times for Gaussian Processes},
  author={Decreusefond, Laurent and Nualart, David},
  journal={The Annals of Probability},
  pages={319--330},
  volume={36},
  number={1},
  year={2008}
}

@article{landau1970supremum,
  title={On the supremum of a Gaussian process},
  author={Landau, Henry J and Shepp, Lawrence A},
  journal={Sankhy{\=a}: The Indian Journal of Statistics, Series A},
  pages={369--378},
  volume={32},
  number={4},
  year={1970}
}

@article{le2008bounded,
  title={On bounded Gaussian processes},
  author={Le, H},
  journal={Statistics \& probability letters},
  volume={78},
  number={6},
  pages={669--674},
  year={2008}
}

@article{bavard1992polygones,
  title={Polygones du plan et polyedres hyperboliques},
  author={Bavard, Christophe and Ghys, {\'E}tienne},
  journal={Geometriae Dedicata},
  volume={43},
  pages={207--224},
  year={1992},
  publisher={Springer}
}

@inproceedings{fillastre2016remark,
author = {Fillastre, Fran{\c{c}}ois and Izmestiev, Ivan},
title = {A remark on spaces of flat metrics with cone singularities of constant sign curvatures},
booktitle = {S{\'e}minaire de th{\'e}orie spectrale et g{\'e}om{\'e}trie},
year = {2016},
pages={65--92}
}

@article{fillastre2017shapes,
  title={Shapes of polyhedra, mixed volumes and hyperbolic geometry},
  author={Fillastre, Fran{\c{c}}ois and Izmestiev, Ivan},
  journal={Mathematika},
  volume={63},
  number={1},
  pages={124--183},
  year={2017},
  publisher={Wiley Online Library}
}

@book{phelps2001lectures,
  title={Lectures on Choquet’s theorem},
  author={Phelps, Robert R},
  year={2001},
  address={New York, NY},
  edition={2},
  publisher={Springer}
}

\noindent{\sc Université Paris-Saclay, Laboratoire de Mathématiques d'Orsay, 91405, Orsay, France}

\noindent{\it Email address:} {\tt \href{mailto:yusen.long@universite-paris-saclay.fr}{yusen.long@universite-paris-saclay.fr}}
\end{document}